\documentclass[10pt,letterpaper]{article}
\usepackage[lmargin=35mm,rmargin=35mm,tmargin=30mm,bmargin=30mm]{geometry}
\usepackage{algorithmic,amsmath,amsthm,graphicx,charter,enumerate}
\usepackage[sort&compress,numbers]{natbib}

\theoremstyle{plain}
\newtheorem{theorem}{Theorem}[section]
\newtheorem{lemma}[theorem]{Lemma}

\newcommand{\twothmref}[2]{Theorems~\ref{thm:#1} and \ref{thm:#2}}
\newcommand{\lemref}[1]{Lemma~\ref{lem:#1}}
\newcommand{\twolemref}[2]{Lemmata~\ref{lem:#1} and \ref{lem:#2}}
\newcommand{\thmref}[1]{Theorem~\ref{thm:#1}}
\newcommand{\figref}[1]{Figure~\ref{fig:#1}}
\newcommand{\secref}[1]{Section~\ref{sec:#1}}

\newcommand{\aaa}{\textup{\hspace*{0.5em}(a)\hspace*{0.25em}}}
\newcommand{\bbb}{\textup{\hspace*{0.5em}(b)\hspace*{0.25em}}}
\newcommand{\ccc}{\textup{\hspace*{0.5em}(c)\hspace*{0.25em}}}

\newcommand{\ONE}{\textup{\hspace*{0.5em}(1)\hspace*{0.25em}}}
\newcommand{\TWO}{\textup{\hspace*{0.5em}(2)\hspace*{0.25em}}}
\newcommand{\THREE}{\textup{\hspace*{0.5em}(3)\hspace*{0.25em}}}

\newcommand{\seclabel}[1]{\label{sec:#1}}
\newcommand{\figlabel}[1]{\label{fig:#1}}

\newcommand{\thmlabel}[1]{\label{thm:#1}}
\newcommand{\lemlabel}[1]{\label{lem:#1}}

\newcommand{\ninth}{\ensuremath{\protect\tfrac{1}{9}}}
\newcommand{\half}{\ensuremath{\protect\tfrac{1}{2}}}
\newcommand{\sixth}{\ensuremath{\protect\tfrac{1}{6}}}
\newcommand{\third}{\ensuremath{\protect\tfrac{1}{3}}}

\newcommand{\ocom}[1]{\ensuremath{\textup{\textsf{ext}}(#1)}}
\newcommand{\icom}[1]{\ensuremath{\textup{\textsf{int}}(#1)}}
\newcommand{\msf}[1]{\ensuremath{\textup{\textsf{msf}}(#1)}}
\newcommand{\F}[2]{\ensuremath{#1\langle#2\rangle}}

\newcommand{\Oh}[1]{\ensuremath{\protect\mathcal{O}(#1)}}

\renewcommand{\baselinestretch}{1.15}
\newcommand{\Figure}[4][!htb]{
\begin{figure}[#1]
	\vspace*{1ex}
	\begin{center}#3\end{center}
	\vspace*{-2ex}
	\caption{#4}
	\figlabel{#2}
\end{figure}
}

\begin{document}

\title{\textbf{Simultaneous Diagonal Flips\\ in Plane Triangulations}\thanks{A preliminary version of this paper was published in the \emph{Proceedings of the 17th Annual ACM-SIAM Symposium on Discrete Algorithms} (SODA '06). Research of the Canadian authors is supported by the Natural Sciences and Engineering Research Council of Canada (NSERC). Research of D.~Wood is supported by the Government of Spain grant MEC SB2003-0270, and by the projects MCYT-FEDER BFM2003-00368 and Gen.\ Cat 2001SGR00224.}}

\renewcommand{\thefootnote}{\fnsymbol{footnote}}

\author{Prosenjit Bose\footnotemark[4]
\and Jurek Czyzowicz\footnotemark[3]
\and Zhicheng Gao\footnotemark[5]
\and Pat Morin\footnotemark[4]
\and David R. Wood\footnotemark[6]}

\footnotetext[4]{School of Computer Science, Carleton University, Ottawa,
Canada (\texttt{\{jit,morin\}@scs.carleton.ca}).}

\footnotetext[3]{D\'{e}partement d'informatique et d'ing\'{e}nierie,
Universit\'{e} du Qu\'{e}bec en Outaouais, Gatineau, Canada (\texttt{Jurek.Czyzowicz@uqo.ca}).}
 
\footnotetext[5]{School of Mathematics and Statistics, Carleton University, 
Ottawa, Canada (\texttt{zgao@math.carleton.ca}).}

\footnotetext[6]{Departament de Matem{\`a}tica Aplicada II, Universitat Polit{\`e}cnica de Catalunya, Barcelona, Spain (\texttt{david.wood@upc.edu}).
Research partially completed at Carleton University and McGill University (Montr\'eal).}

\maketitle

\begin{abstract}  
Simultaneous diagonal flips in plane triangulations are investigated. It is proved that every $n$-vertex triangulation with at least six vertices has a simultaneous flip into a $4$-connected triangulation, and that
it can be computed in \Oh{n} time. It follows that every triangulation has a simultaneous flip into a Hamiltonian triangulation. This result is used to prove that for any two $n$-vertex triangulations, there exists a sequence of \Oh{\log n} simultaneous flips to transform one into the other. The total number of  edges flipped in this sequence is \Oh{n}. The maximum size of a simultaneous flip is then studied. It is proved that every triangulation has a simultaneous flip of at least $\tfrac{1}{3}(n-2)$ edges. On the other hand, every simultaneous flip has at most $n-2$ edges, and there exist triangulations with a maximum simultaneous flip of $\tfrac{6}{7}(n-2)$ edges.\\[1ex] \textbf{keywords:} graph, plane triangulation, diagonal flip, simultaneous flip, Hamiltonian 
\end{abstract}

%%%%%%%%%%%%%%%%%%%%%%%%%%%%%%%%%%%%%%%%%%%%%%%%%%%%%%%%%%%%%%%%%%%%%%%%%%%%%
\section{Introduction}
\seclabel{Introduction}
%%%%%%%%%%%%%%%%%%%%%%%%%%%%%%%%%%%%%%%%%%%%%%%%%%%%%%%%%%%%%%%%%%%%%%%%%%%%%

A (\emph{plane}) \emph{triangulation} is a simple planar graph with a fixed (combinatorial) plane embedding in which every face is bounded by a \emph{triangle} (that is, a $3$-cycle). So that we can speak of the interior and exterior of a cycle, one face is nominated to be the \emph{outerface}, although the choice of outerface is not important for our results. 

Let $vw$ be an edge of a triangulation $G$. Let $(v,w,x)$ and $(w,v,y)$ be the faces incident to $vw$. Then $x$ and $y$ are distinct vertices, unless $G=K_3$. We say that $x$ and $y$ \emph{see} $vw$. Let $G'$ be the embedded graph obtained from $G$ by deleting $vw$ and adding the edge $xy$, such that in the cyclic order of the edges incident to $x$ (respectively, $y$), $xy$ is added between $xv$ and $xw$ ($yw$ and $yv$). If $G'$ is a triangulation, then $vw$ is (\emph{individually}) \emph{flippable}, and $G$ is \emph{flipped} into $G'$ by $vw$. This operation is called a (\emph{diagonal}) \emph{flip}, and is illustrated in \figref{ExampleSimultaneousFlip}. If $G'$ is not a triangulation and $G\ne K_3$, then $xy$ is already an edge of $G$; we say that $vw$ is \emph{blocked} by $xy$, and $xy$ is a \emph{blocking} edge. 

\Figure{ExampleSimultaneousFlip}{\includegraphics{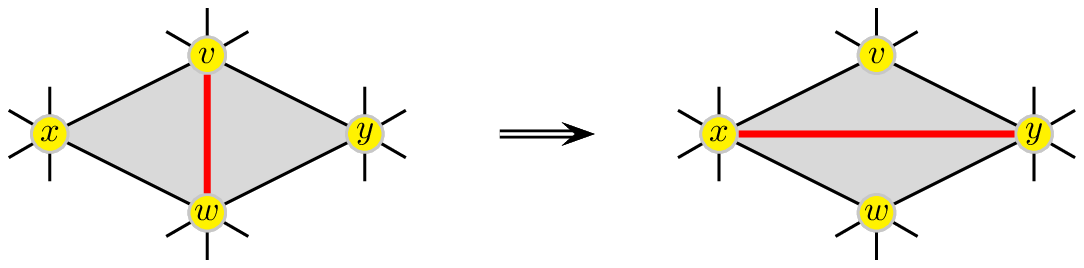}}{Edge
$vw$ is flipped into $xy$.}

In 1936, \citet{Wagner36} proved that a finite sequence of diagonal flips transform a given triangulation into any other triangulation with the same number of vertices. Since then diagonal flips in plane triangulations \citep{Komuro-Yoko97, NSS-JGT06, Negami-DM94, STT-SJDM92, GUW-GC01, KA-AC01, Negami99, Nakamigawa-TCS00, HOS-JUCS96, MNO-GC03, GW-JCTA99} and in triangulations of other surfaces \citep{CGMN-DM02, Negami99, CN-DM00a, WN-Yoko99, KNN-JCTB99, Negami-Yoko98, BNN-JCTB96, Negami-DM94, NN93, NN-Yoko02} have been studied extensively. It can be shown that for triangulation with $n$ vertices, the number of flips in Wagner's proof is \Oh{n^2}. \citet{Komuro-Yoko97} improved this bound to \Oh{n}. The best known bound is $\max\{6n-30,0\}$ due to \citet{MNO-GC03}. 

For labelled triangulations, \citet{STT-SJDM92} proved that \Oh{n \log n} flips are sufficient to transform one labelled triangulation with $n$ vertices into any other, and $\Omega(n \log n)$ flips are sometimes necessary. The upper bound was independently rediscovered by \citet{GUW-GC01}. Note that the above-mentioned \Oh{n} upper bound in the unlabelled setting \citep{Komuro-Yoko97, MNO-GC03} can also be obtained by a careful analysis of the proof by \citet{STT-SJDM92}. 

\citet{Wagner36} in fact proved that every $n$-vertex triangulation can be transformed by a sequence of flips into the so-called \emph{standard} triangulation $\Delta_n$, which is illustrated in \figref{StandardTriangulation} and is defined as the triangulation on $n$ vertices with two \emph{dominant} vertices (adjacent to every other vertex). Clearly two $n$-vertex triangulations each with two dominant vertices are isomorphic. To transform one $n$-vertex triangulation $G_1$ into another $G_2$, first transform $G_1$ into $\Delta_n$, and then apply the flips to transform $G_2$ into $\Delta_n$ in reverse order. A similar approach is used in this paper in the context of simultaneous flips in triangulations. 

\Figure{StandardTriangulation}{\includegraphics{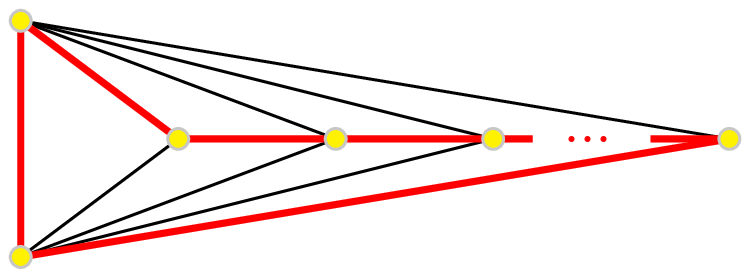}}{The
standard triangulation and a Hamiltonian cycle.}

Let $S$ be a set of edges in a plane triangulation $G$. The embedded graph obtained from $G$ by flipping every edge in $S$ is denoted by \F{G}{S}. If \F{G}{S} is a triangulation, then $S$ is (\emph{simultaneously}) \emph{flippable} in $G$, and $G$ is \emph{flipped} into \F{G}{S} by $S$. This operation is called a \emph{simultaneous} (\emph{diagonal}) \emph{flip}. Note that it is possible for $S$ to be flippable, yet $S$ contains non-flippable edges, and it is possible for every edge in $S$ to be flippable, yet $S$ itself is not flippable. As far as the authors are aware, simultaneous flips have previously been studied only in the more restrictive context of geometric triangulations of a point set \citep{GHNP-IJCGA03}. Individual flips have also been studied in a geometric context \citep{HNU-DCG99, HN-CGTA99}.

In \secref{Basics} we characterise flippable sets and give a number of introductory lemmas. Our first main result states that every triangulation with at least six vertices can be transformed by one simultaneous flip into a $4$-connected (and hence Hamiltonian) triangulation. Moreover, this flip can be computed in \Oh{n} time for $n$-vertex triangulations. These results are presented in \secref{Hamiltonian}. In \secref{Outerplane} we study simultaneous flips in maximal outerplanar graphs. We prove that for any two $n$-vertex maximal outerplanar graphs, there exists a sequence of \Oh{\log n} simultaneous  flips to transform one into the other. The method used is the basis for the main result in \secref{TwoTriangulations}, which states that for any two $n$-vertex triangulations, there exists a sequence of \Oh{\log n} simultaneous  flips to transform one into the other. This result is optimal for many pairs of triangulations. For example, if one triangulation has $\Theta(n)$ maximum degree and the other has \Oh{1} maximum degree, then $\Omega(\log n)$ simultaneous  flips are needed, since one simultaneous  flip can at most halve the degree of a vertex. This also holds for diameter instead of maximum degree. Finally in \secref{BigFlip} the maximum size of a simultaneous flip is studied. It is proved that every triangulation has a simultaneous flip of at least $\third(n-2)$ edges. On the other hand, every simultaneous flip has at most $n-2$ edges, and there exist triangulations with a maximum simultaneous flip of $\tfrac{6}{7}(n-2)$ edges.

%%%%%%%%%%%%%%%%%%%%%%%%%%%%%%%%%%%%%%%%%%%%%%%%%%%%%%%%%%%%%%%%%%%%%%%%%%%%%
\section{Basics}
\seclabel{Basics}
%%%%%%%%%%%%%%%%%%%%%%%%%%%%%%%%%%%%%%%%%%%%%%%%%%%%%%%%%%%%%%%%%%%%%%%%%%%%%

We start with a characterisation of flippable sets that is used throughout the paper. Two edges of a triangulation that are incident to a common face are \emph{consecutive}. If two consecutive edges are simultaneously flipped, then the two new edges cross, as illustrated in \figref{Obstacles}(a). Thus no two edges in a flippable set are consecutive. Two edges form a \emph{bad pair} if they are seen by the same pair of vertices. If a bad pair of edges are simultaneously flipped, then the two new edges are parallel, as illustrated in \figref{Obstacles}(b). Thus no two edges in a flippable set form a bad pair. If an edge $vw$ is blocked by an edge $pq$ as illustrated in \figref{Obstacles}(c), then $vw$ is not individually flippable, but $vw$ can be in a flippable set $S$ as long as $pq$ is also in $S$. We now show that these three properties characterise flippable sets.

\Figure{Obstacles}{\includegraphics{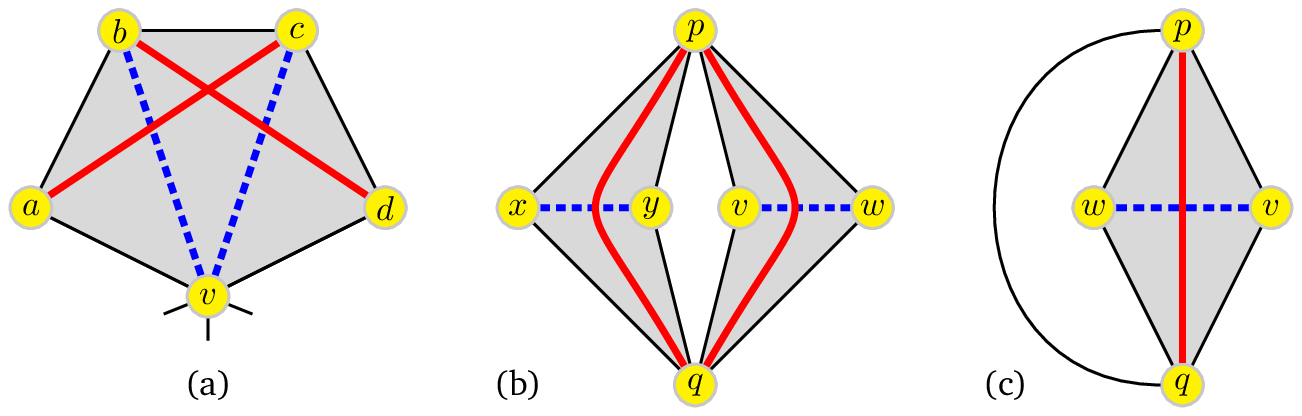}}{Obstacles to a flippable set. Dashed edges are flipped to create bold edges. Shaded regions are faces}

\begin{lemma}
\lemlabel{FlippableCharacterisation}
A set of edges $S$ in a triangulation $G\ne K_3$ is flippable if and only if:\\
\ONE\ no two edges in $S$ are consecutive, \\
\TWO\ no two edges in $S$ form a bad pair, and\\
\THREE\ for every edge $vw\in S$, either $vw$ is flippable or the edge that blocks $vw$ is also in $S$.
\end{lemma}

\begin{proof}
We have already seen that each condition is necessary for $S$ to be flippable. Now suppose that all three conditions are satisfied. Since no two edges in $S$ are consecutive, \F{G}{S} is a graph embedded in the plane and every face is a triangle. Suppose that two edges $e_1$ and $e_2$ are parallel in \F{G}{S}. Since $G$ has no parallel edges, $e_1$ and $e_2$ are both not in $G$. If exactly one of $e_1$ and $e_2$ is in $G$, then condition (3) fails. If neither of $e_1$ and $e_2$ are in $G$, then the edges in $S$ that flipped to $e_1$ and $e_2$ form a bad pair.
\end{proof}

Note that condition (1) in \lemref{FlippableCharacterisation} is equivalent to saying that the edges of the dual $G^*$ that correspond to $S$ form a matching. (The \emph{dual} $G^*$ of $G$ is the plane graph with one vertex for every face of $G$, such that two vertices of $G^*$ are adjacent whenever the corresponding faces in $G$ are incident to a common edge.)\ 

A cycle $C$ in a triangulation $G$ is \emph{separating} if deleting the vertices of $C$ from $G$ produces a disconnected graph. 

%%%%%%%%%%%%%%%%%%%%%%%%%%%%%%%%%%%%

\begin{lemma}
\lemlabel{EdgeInSepTriangle}
An edge in a separating triangle $T$ of a triangulation is individually flippable.
\end{lemma}

\begin{proof}
Consider an edge $vw$ in $T$. Say $vw$ is seen by $p$ and $q$. Then one of $p$ and $q$ is inside $T$, and the other is outside $T$. Thus $pq$ is not an edge, and $vw$ is flippable. 
\end{proof}

%%%%%%%%%%%%%%%%%%%%%%%%%%%%%%%%%%%

The next observation quickly follows from the Jordan Curve Theorem.

\begin{lemma}
\lemlabel{EdgeInCycle}
Let $vw$ be an edge of a triangulation that is seen by vertices $p$ and $q$. Suppose that $p$ is inside some cycle $C$ and $q$ is outside $C$. Then $vw\in C$.\qed
\end{lemma}

%%%%%%%%%%%%%%%%%%%%%%%%%%%%%%%%%

The next two results show that blocking edges are nearly always flippable, and except for essentially one case, do not appear in a bad pair.

\begin{lemma}
\lemlabel{BlockingEdgeFlippable}
A blocking edge is individually flippable in a triangulation $G\ne K_4$.
\end{lemma}

\begin{proof} 
Let $vw$ be an edge of $G$ that is blocked by $pq$. Without loss of generality, $w$ is inside the triangle $pvq$. If $pvq$ is a separating triangle, then $pq$ is flippable by \lemref{EdgeInSepTriangle}. If $pvq$ is not separating, then $pwq$ must be a separating triangle since $G\ne K_4$.  Therefore, $pq$ is flippable by \lemref{EdgeInSepTriangle}.
\end{proof}

%%%%%%%%%%%%%%%%%%%%%%%%%%%%%%

\begin{lemma}
\lemlabel{BlockingEdgeInBadPair}
Suppose that $vw$ and $xy$ are a bad pair in a triangulation $G$, both seen by vertices $p$ and $q$. Suppose that $vw$ blocks some edge $ab$. Then $xy$ and $ab$ are consecutive, and $vw$ and $xy$ are in a common triangle (amongst other properties). 
\end{lemma}

\begin{proof} Without loss of generality, $w$ and $x$ are inside the cycle $(v,p,y,q)$, and $b$ is inside the triangle $(v,a,w)$, as illustrated in \figref{BlockingInBadPair}. Now $(v,p,y,q)$ is a separating $4$-cycle with $w$ in its interior. Since $w$ is adjacent to $a$ and $b$, both $a$ and $b$ must be on the boundary of $(v,p,y,q)$. It follows that $p=b$ and $y=a$. If $w\ne x$, then the neighbours $w$ and $a$ are respectively on the inside and outside of the cycle $(p,x,q,v)$, which is not possible. Thus $w=x$. Hence $xy$ and $ab$ are consecutive, and $vw$ and $xy$ are in a common triangle. \end{proof}

\Figure{BlockingInBadPair}{\includegraphics{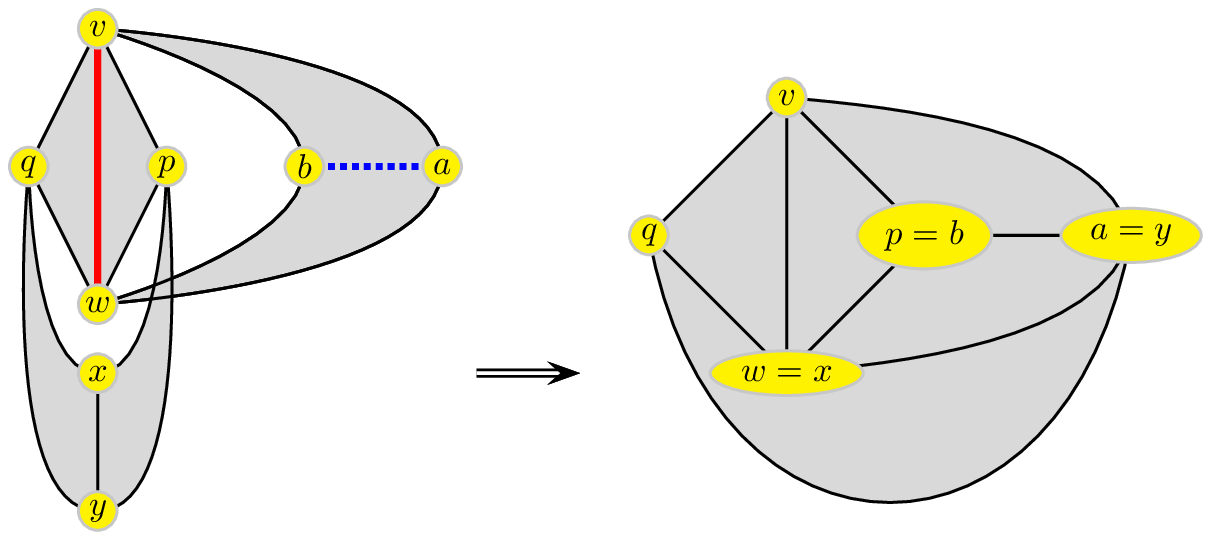}}{The only case when a blocking edge $vw$ is in a bad pair.}

%%%%%%%%%%%%%%%%%%%%%%%%%%%%%%%%%%%%%%%%%%%%%%%%%%%%%%%%%%%%%%%%%%%%%%%%%%%%%
\section{Flipping into a $4$-Connected Triangulation}
\seclabel{Hamiltonian}
%%%%%%%%%%%%%%%%%%%%%%%%%%%%%%%%%%%%%%%%%%%%%%%%%%%%%%%%%%%%%%%%%%%%%%%%%%%%%

The main result in this section is that every triangulation with at least six vertices has a simultaneous flip into a $4$-connected (and hence Hamiltonian) triangulation. It is well known that a triangulation is $4$-connected if and only if it has no separating triangle. Thus our focus is on flips that break separating cycles.

\begin{lemma}
\lemlabel{BasicLemma}
Let $S$ be a set of edges in a triangulation such that no two edges of $S$ are in a common triangle, and every edge in $S$ is in a separating triangle.  Then $S$ is flippable.
\end{lemma}

\begin{proof} 
By \lemref{EdgeInSepTriangle}, every edge in $S$ is individually flippable. 
Thus, by \lemref{FlippableCharacterisation}, it suffices to prove that no two edges in $S$ form a bad pair. Suppose that $vw,xy\in S$ form a bad pair. Then $vw$ and $xy$ are seen by the same pair of vertices $p$ and $q$. Let $T$ be a separating triangle containing $vw$. Then one of $p$ and $q$ is inside $T$, and the other is outside $T$. By \lemref{EdgeInCycle}, $xy$ must be an edge of $T$, which implies that $vw$ and $xy$ are in a common triangle. This contradiction proves that there is no bad pair of edges both in $S$, and $S$ is flippable.
%
%Since no two edges of $S$ appear in the same triangle, no two edges in $S$ are consecutive. Suppose on the contrary that \F{G}{S} has parallel edges $e_1=vw$ and $e_2=vw$. Let $S'$ be the set of edges in $\F{G}{S}$ that are not in $G$. 
%
%Initially suppose that exactly one of $e_1$ and $e_2$, say $e_1$, is in $S'$; see \figref{NoParallelEdges}(a). Let $xy$ be the edge of $G$ flipped to $e_1$. By assumption, $xy$ is in a separating triangle $xyz$ of $G$. Then in $G$, $z$ is either inside the triangle $\{e_2,vx,xw\}$ or outside the triangle $\{vy,yw,e_2\}$. ($z\ne v$ and $z\ne w$ as otherwise $xyz$ would not be a separating triangle.)\ However, in these cases, $z$ could not be adjacent to $y$ and $x$, respectively, without breaking planarity.  
%
%Now suppose that $e_1$ and $e_2$ are both in $S'$; see \figref{NoParallelEdges}(b). Let $xy$ be the edge of $G$ flipped to $e_1$, and let $rs$ be the edge of $G$ flipped to $e_2$. By assumption, $rs$ is in a separating triangle $rst$ of $G$. Thus $rt,st\not\in S$. Hence $v$ is on the inside of $rst$ and $w$ is on the outside of $rst$ (without loss of generality). Thus a vertex adjacent to both $v$ and $w$ is in the separating triangle $rst$. Since $x$ and $y$ are both adjacent to $v$ and $w$, the edge $xy$ is in the separating triangle $rst$. Thus $xy=rs$, which is a contradiction.
%Hence $S$ is flippable. 
\end{proof}

%\Figure{BasicLemmaFigure}{\includegraphics{BasicLemmaFigure}}{Illustration for \lemref{BasicLemma}.}

%By assumption, each of $xy$ and $rs$ is in a separating triangle of $G$. Let $z$ and $t$ be the vertices of $G$ such that $xyz$ and $rst$ are separating triangles. For $xyz$ and $rst$ to be separating triangles of $G$, it must be the case that $z=r$, $t=x$, and $y=s$; or $z=s$, $t=y$, and $x=r$. In either case, $xy$ and $rs$ are in a common triangle in $G$, contradicting the initial assumption. In each case, we have a contradiction. 

%\Figure{NoParallelEdges}{\includegraphics{NoParallelEdges}}{Dashed edges are flipped to create bold parallel edges. }

%\lemref{BasicLemma} gives a sufficient condition for a set of edges $S$
%to be flippable; namely, every edge in $S$ is in a separating triangle. This
%condition can be dropped if $G$ is $5$-connected.

\begin{lemma}
\lemlabel{IsFourConnected}
Let $G$ be a triangulation with $n\geq6$ vertices. Let $S$ be a set of edges in $G$ that satisfy the conditions in \lemref{BasicLemma}, and suppose that every separating triangle contains an edge in $S$. Then \F{G}{S} is $4$-connected.
\end{lemma}

\begin{proof} Suppose on the contrary, that \F{G}{S} contains a
separating triangle $T=(u,v,w)$. Let $S'$ be the set of edges in $\F{G}{S}$ that
are not in $G$. We proceed by case-analysis on $|T\cap S'|$ (refer to
\figref{NoSeparatingTriangle}). Since every separating triangle in $G$ has an
edge in $S$, $|T\cap S'|\geq1$.

Case 1. $|T\cap S'|=1$: Without loss of generality, $vw\in S'$, $uv\not\in S'$, and $uw\not\in S'$. Suppose $xy$ was flipped to $vw$. Then $xy$ is in a separating triangle $xyp$ in $G$. Any vertex adjacent to both $v$ and $w$ must be a vertex of the separating triangle $xyp$. Thus $p=u$. Since $G$ has at least six vertices, at least one of the triangles $\{(u,v,x),(u,v,y),(u,w,x),(u,w,y)\}$ is a separating triangle. Thus at least one of the edges in these triangles is in $S$.  Since $xy\in S$, and no two edges of $S$ appear in a common triangle, $\{ux,uy,vx,vy,wx,wy\}\cap S=\emptyset$. Thus $uv$ or $uw$ is in $S$. But then $uvw$ is not a triangle in \F{G}{S}, which is a contradiction.

\Figure{NoSeparatingTriangle}{\includegraphics{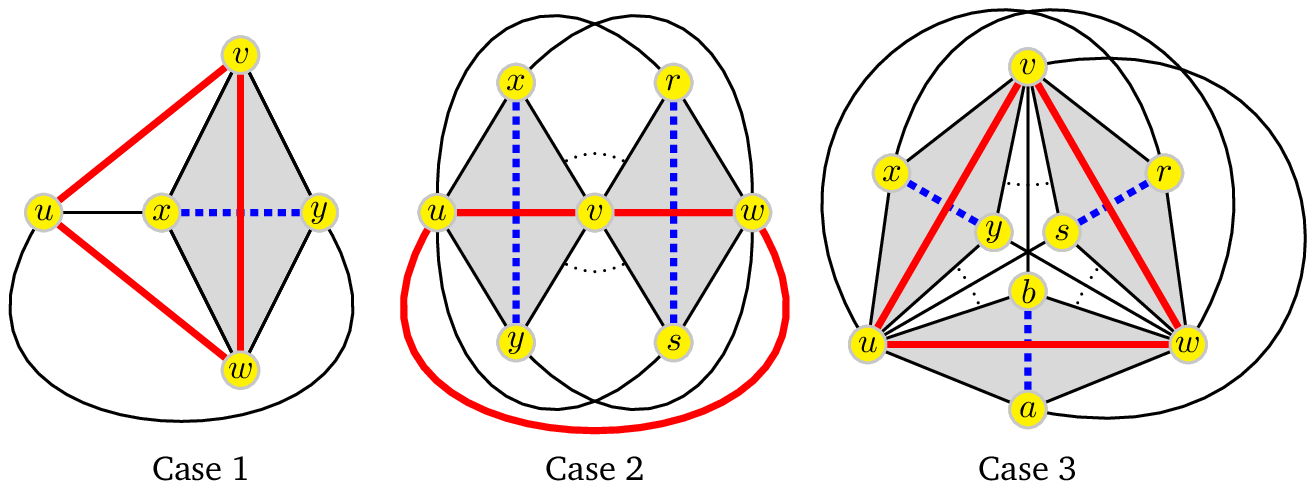}}{Dashed edges are flipped to create a bold separating triangle. Shaded regions are faces.}

Case 2. $|T\cap S'|=2$: Without loss of generality, $uv\in S'$, $vw\in S'$, and $uw\not\in S'$. Suppose $xy$ was flipped to $uv$, and $rs$ was flipped to $vw$. Without loss of generality, $y$ and $s$ are inside $uvw$ in \F{G}{S}. Then in $G$, $xy$ was in a separating triangle $xyz$, and $rs$ was in a separating triangle $rst$.  By an argument similar to that in \lemref{BasicLemma},  $z=w$ and $t=u$. But then the subgraph of $G$ induced by $\{u,v,w,x,y,r,s\}$ is not planar, or it contains parallel edges in the case that $x=r$ and $y=s$.

Case 3. $|T\cap S'|=3$: Suppose $xy$ was flipped to $uv$, $rs$ was flipped to $vw$, and  $ab$ was flipped to $uw$. Without loss of generality, $y$, $s$ and $b$ are inside $uvw$ in \F{G}{S}.  In $G$, $xy$ was in a separating triangle $xyz$, $rs$ was in a separating triangle $rst$, and $ab$ was in a separating triangle $(a,b,c)$. By an argument similar to that in \lemref{BasicLemma}, $z=w$, $t=u$, and $c=v$. But then the subgraph of $G$ induced by $\{u,v,w,x,y,r,s,a,b\}$ is non-planar, or contains parallel edges in the case
that $y=s=b$ and $x=r=a$.

In each case we have derived a contradiction. Therefore \F{G}{S} has no separating triangle, and thus is $4$-connected.
\end{proof}

%%%%%%%%%%%%%%%%%%%%%%%%%%%%%%%%%%%%%%%%%%%%%%%%%%%%%%%%%%%%%%%%%%%%%%%%%%%%%%%

Observe that the restriction in \lemref{IsFourConnected} to triangulations with
at least six vertices is unavoidable. Every triangulation with at most five
vertices has a vertex of degree three, and is thus not $4$-connected. 

We now consider how to determine a set of edges that satisfy \lemref{IsFourConnected}. 

\begin{lemma}
\lemlabel{MatchingInPrimal}
Let $e$ be an edge of an $n$-vertex triangulation $G$. Then $G$ has a set of
edges $S$ that can be computed in \Oh{n} time such that $e\in S$ and every face
of $G$ has exactly one edge in $S$.
\end{lemma}

\begin{proof}
\citet{BBDL-JAlg01} proved the following strengthening of Petersen's matching theorem \citep{Petersen1891}: every $3$-regular bridgeless planar graph has a perfect matching that contains a prespecified edge and can be computed in linear time. The dual $G^*$ is a $3$-regular bridgeless planar graph with $2n-4$ vertices. A perfect matching in $G^*$ corresponds to the desired set $S$.
\end{proof}

\lemref{MatchingInPrimal} only accounts for triangles of $G$ that are faces. We account for separating triangles as follows.

\begin{lemma}
\lemlabel{TriangleCover}
Let $e$ be an edge of an $n$-vertex triangulation $G$. Then $G$ has a set of
edges $S$ such that $e\in S$ and every triangle of $G$ has exactly one edge in
$S$.
\end{lemma}

\begin{proof}  
We proceed by induction on the number of separating triangles. The result follows for a triangulation with no separating triangles by \lemref{MatchingInPrimal}. Now suppose $G$ has $k>0$ separating triangles, and the lemma holds for triangulations with less than $k$ separating triangles. Let $T$ be a separating triangle of $G$. Let the components of $G\setminus T$ have vertex sets $V_1$ and $V_2$. Consider the induced subgraphs $G_1:=G[V_1\cup T]$ and $G_2:=G[V_2\cup T]$. Without loss of generality, the given edge $e$ is in $G_1$.  Both $G_1$ and $G_2$ have less than $k$ separating triangles. By induction $G_1$ has a set of edges $S_1$ such that $e\in S_1$, and every triangle of $G_1$ has exactly one edge in $S_1$. Let $e_2$ be the edge in $S_1\cap T$. By induction, $G_2$ has a set of edges $S_2$ such that $e_2\in S_2$, and every triangle of $G_2$ has exactly one edge in $S_2$.  Thus $S:=S_1\cup S_2$ is a set of edges of $G$ such that $e\in S$, and every triangle of $G$ has exactly one edge in $S$. 
\end{proof}

By taking as a flip set those edges in the set $S$ from \lemref{TriangleCover} that are in some separating triangle, \lemref{IsFourConnected} implies that every triangulation with at least six vertices  has a simultaneous flip into a $4$-connected triangulation. However, due to the presence of separating triangles, it is not obvious how to implement  \lemref{TriangleCover} in linear  time. In what follows we show how to do this.

First we outline a few properties of separating triangles.  Let $T$ be a separating triangle of a triangulation $G$. Removing the vertices of $T$  from $G$ produces two components, an \emph{inner component} (containing no vertex on the outerface) and an \emph{outer component}. Denote by \icom{T} and \ocom{T} the sets of vertices of the inner and outer components. Define a \emph{containment} relation, denoted by $\preceq$, on the set of separating triangles of $G$ as follows. For two separating triangles $T_1$ and $T_2$ of $G$, let $T_1 \preceq T_2$ whenever $\icom{T_1} \subseteq  \icom{T_2}$. Clearly $\preceq$ is a  partial order. 

We first show how to compute a linear extension $R$ of $\preceq$ in linear time. We then show how to use $R$ to compute the set $S$ in \lemref{TriangleCover} in linear time. The \emph{canonical ordering} of \citet{dFPP90} will be a useful tool. Let $G$ be a plane triangulation with outerface $(a,b,c)$. A linear ordering of the vertices $(v_1=a,v_2=b,v_3,\dots,v_n=c)$ is \emph{canonical} if the following conditions hold for all $3\leq i\leq n$:

\begin{itemize}
\item the subgraph $G_i$ induced by $\{v_1,v_2,\dots,v_i\}$ is $2$-connected, and the boundary of its outerface is a cycle $C_i$ containing the edge $ab$; and
\item the vertex $v_i$ is in the outerface of $G_{i-1}$, and the  neighbours of $v_i$ in $G_{i-1}$ form a subinterval of the path $C_{i-1}\setminus\{ab\}$ consisting of at least two vertices (and $v_3$ is adjacent to $v_1$ and $v_2$).
\end{itemize} 

\citet{dFPP90} proved that every triangulation has a canonical ordering. Define the \emph{level} of a separating triangle $T$, denoted by $\ell(T)$, as the largest index of a vertex of $T$ in a given canonical order.

\begin{lemma} 
\lemlabel{level}
Let $T_1$ and $T_2$ be separating triangles  such that $\ell(T_1)<\ell(T_2)$.
Then $T_1 \preceq T_2$ or  $\icom{T_1}\cap\icom{T_2}=\emptyset$.  
\end{lemma}

\begin{proof} 
Let $T_1=(a,b,c)$ and $T_2=(x,y,z)$. Suppose on the contrary that $T_2 \preceq T_1$. Then $\icom{T_2}\subset\icom{T_1}$ since $T_1$ and $T_2$ are distinct. Without loss of generality, let $c$ be the vertex of $T_1$ defining $\ell(T_1)=i$, and let $z$ be the vertex of $T_2$ defining $\ell(T_2)=j$. Since $\ell(T_1)=i<\ell(T_2)=j$, $c$ is distinct from $z$. By the canonical ordering, no vertex in $\icom{T_1}$ is on the outerface of any $G_k$ for $k \geq i$. Since $z$ occurs after $c$ in the canonical ordering, all the vertices adjacent to $z$ in $\icom{T_2}$ are on the outerface of $G_i$. This implies that none of these vertices are in $\icom{T_1}$, which is the desired contradiction.  
\end{proof}

\begin{lemma}
\lemlabel{totalorder}
For an $n$-vertex plane triangulation $G$, a linear extension $R$ of $\preceq$ 
can be computed in \Oh{n} time.
\end{lemma}

\begin{proof}  
First note that a canonical ordering can be computed in \Oh{n} time \citep{dFPP90} (also see \citep{ChrobakPayne-IPL95}). \lemref{level} implies that if all of the separating triangles of $G$ have different levels, then ordering them by increasing level gives the desired linear extension $R$. What remains is to order the separating triangles at the same level. These triangles share a common vertex $v_i$ that defines their level. The neighbours of $v_i$ in $G_{i-1}$ form a path $P=(p_1,p_2,\ldots,p_k)$ on the boundary of the outerface of $G_{i-1}$.  Every separating triangle of $G$ at level $i$ consists of $v_i$ and two non-consecutive vertices of $P$. To establish the containment relation between these triangles, we simply need to look at the indices of the vertices of $P$. Let $T_1=(v_i,p_a,p_b)$ and $T_2=(v_i,p_c,p_d)$ be distinct separating triangles with $a<b$ and $c<d$.  If $a<b\leq c<d$ or $c<d\leq a<b$  then $\icom{T_1}\cap\icom{T_2}=\emptyset$ by the canonical ordering. It is impossible for $a<c<b<d$ or $c<a<d<b$ since the graph induced on $P$ is outerplanar and this would violate planarity. If $a \leq c < d \leq b$ then $T_2 \preceq T_1$ and if $c \leq a < b \leq d$ then $T_1 \preceq T_2$.  Since we can compute the graph induced by $\{p_1, p_2, \ldots, p_k\}$ in \Oh{k} time, all of the separating triangles at level $i$ can be ordered in \Oh{k} time by performing a breadth-first search on the graph induced on $P$. The result follows since the sum of the degrees of a plane graph is \Oh{n}. 
\end{proof}

We now turn our attention to computing the set $S$ from \lemref{TriangleCover} in linear time. Denote by \textsc{FaceSet}($G,e$) the set $S$ from \lemref{MatchingInPrimal}; that is, every face of $G$ has exactly one edge in $S$, and if $e$ is specified then $e\in S$. 

\medskip\noindent\textbf{Algorithm} \textsc{TriangleSet}$(G,R)$

\noindent\emph{Input:} triangulation $G$, ordered list of separating triangles $R$ of $G$.

\noindent\emph{Output:} a set $S$ of edges of $G$ such that every triangle of $G$ has exactly one edge in $S$.

\begin{algorithmic}[1]

\IF{$R=\emptyset$}\STATE{return \textsc{FaceSet}($G$, unspecified);}
\ELSE
\STATE{let $T$ be the first triangle in $R$;}
\STATE{let $S:=$\textsc{TriangleSet}($G\setminus \icom{T}, R\setminus T$);}
\STATE{let $e$ be the edge in $S\cap T$;}
\STATE{return $S\cup\textsc{FaceSet}(G\setminus \ocom{T}, e)$;}
\ENDIF

\end{algorithmic}

We now prove the correctness and running time of the algorithm.

\begin{lemma}
\lemlabel{AlgorithmicTriangleCover}
For every $n$-vertex triangulation $G$, the algorithm \textsc{TriangleSet}$(G,R)$ returns a set $S$ consisting of exactly one edge in every triangle of $G$. The running time is \Oh{n}.
\end{lemma}

\begin{proof} We proceed by induction on $|R|$. If $R=\emptyset$ then every
triangle in $G$ is a face, and \textsc{TriangleSet}$(G,R)$ correctly
computes $S$ with a call to \textsc{FaceSet}$(G)$. Now assume that
$R\ne\emptyset$. Let $T$ be the first triangle in $R$. Then $T$ is an innermost
separating triangle of $G$, and $G\setminus\ocom{T}$ has no separating triangle.
Hence $R\setminus T$ is a linear extension of the containment relation $\preceq$ on the set of separating triangles of $G\setminus\icom{T}$. By induction, $S:=\textsc{TriangleSet}(G\setminus \icom{T},R\setminus T)$ consists of exactly one edge in every triangle of $G\setminus\icom{T}$. Thus there is exactly one edge $e\in S\cap T$. Every triangle in $G\setminus\ocom{T}$ is a face. By \lemref{MatchingInPrimal}, $\textsc{FaceSet}(G\setminus\ocom{T},e)$ consists of exactly one edge in every triangle of $G\setminus\ocom{T}$ including $e$. Together with $S$ we have the desired set for $G$. The running time is described by the recurrence $X(n)=X(n-|\icom{T}|)+\Oh{|\icom{T}|}+\Oh{1}$, which solves to \Oh{n}.
\end{proof}

Note that Algorithm \textsc{TriangleSet} can be easily modified to guarantee that a prespecified edge is in $S$. 

\begin{theorem} 
\thmlabel{SingleFlip}
Every triangulation $G$ with $n\geq6$ vertices  has a simultaneous  flip into a
$4$-connected triangulation that can be computed in \Oh{n} time.
\end{theorem}

\begin{proof} 
By \lemref{AlgorithmicTriangleCover}, $G$ has a set of edges $S$
such that every separating triangle of $G$ has exactly one edge in $S$ and no
triangle of $G$ contains two edges of $S$. By \lemref{BasicLemma}, $S$ is
flippable. By \lemref{IsFourConnected}, \F{G}{S} is $4$-connected. 
\end{proof}

%%%%%%%%%%%%%%%%%%%%%%%%%%%%%%%%%%%%%%%%%%%%%%%%%%%%%%%%%%%%%%%%%%%%%%%%%%%%

We can obtain a stronger result at the expense of a slower algorithm. The
following consequence of the $4$-colour theorem is essentially a Tait
edge-colouring \citep{Tait1880a}.

\begin{lemma} 
\lemlabel{TaitColouring}
Every $n$-vertex planar graph $G$ has an edge $3$-colouring that can be
computed in \Oh{n^2} time, such that every triangle is trichromatic.
\end{lemma}

\begin{proof} \citet{RSST97} proved that $G$ has a proper vertex 4-colouring
that can be computed in \Oh{n^2} time. Let the colours be $\{1,2,3,4\}$.
Colour an edge \emph{red} if its endpoints are coloured $1$ and $2$, or $3$
and $4$. Colour an edge \emph{blue} if its endpoints are coloured $1$ and
$3$, or $2$ and $4$. Colour an edge \emph{green} if its endpoints are
coloured $1$ and $4$, or $2$ and $3$. Since the vertices of each triangle $T$ are trichromatic, the edges of $T$ are also trichromatic.
\end{proof}

%%%%%%%%%%%%%%%%%%%%%%%%%%%%%%%%%%%%%%%%%%%%%%%%%%%%%%%%%%%%%%%%%%%%%%%%%%%%%%

\begin{theorem}   
\thmlabel{ThreeFlips}  
Let $G$ be a triangulation with $n\geq6$ vertices. Then $G$ has three pairwise
disjoint flippable sets of edges $S_1,S_2,S_3$ that can be computed in \Oh{n^2}
time, such that each \F{G}{S_i} is $4$-connected.
\end{theorem}

\begin{proof}  By \lemref{TaitColouring}, $G$ has an edge $3$-colouring such that every triangle is trichromatic. For any of the three colours, let $S$ be the set of edges receiving that colour and in a separating triangle. By \lemref{BasicLemma}, $S$ is flippable. By \lemref{IsFourConnected}, \F{G}{S} is $4$-connected.  \end{proof}

We have the following corollary of \twothmref{SingleFlip}{ThreeFlips}, since
every triangulation on at most five vertices (that is, $K_3$, $K_4$ or
$K_5\setminus e$) is Hamiltonian, and every $4$-connected triangulation has a
Hamiltonian cycle \citep{Whitney31} that can be computed in linear time
\citep{CN89}.

\begin{theorem}
\thmlabel{Hamiltonian} 
Every $n$-vertex triangulation $G$ has a simultaneous  flip into
a Hamiltonian triangulation that can be computed in \Oh{n} time. Furthermore,
$G$ has three disjoint simultaneous flips  that can be computed in \Oh{n^2} time, such that each transforms $G$ into a Hamiltonian triangulation. \qed
\end{theorem}

%%%%%%%%%%%%%%%%%%%%%%%%%%%%%%%%%%%%%%%%%%%%%%%%%%%%%%%%%%%%%%%%%%%%%%%%%%%%%
\section{Outerplane Graphs}
\seclabel{Outerplane}
%%%%%%%%%%%%%%%%%%%%%%%%%%%%%%%%%%%%%%%%%%%%%%%%%%%%%%%%%%%%%%%%%%%%%%%%%%%%%

A plane graph is \emph{outerplane} if every vertex lies on the outerface. The other faces are \emph{internal}. An edge that is not on the boundary of the outerface is \emph{internal}. Let $G$ be an (edge-)maximal outerplane graph $G$ on $n$ vertices. Every internal face is a triangle, and $G$ has $2n-3$ edges and $n-2$ internal faces. The \emph{dual tree} of $G$, denoted by $G^*$, is the dual graph of $G$ without a vertex corresponding to the outerface. Observe that $G^*$ is a tree with $n-2$ vertices and maximum degree at most three.

The notions of \emph{diagonal flip} and \emph{flippable set} for triangulations
are extended to maximal outerplane graphs in the natural way, except that only
internal edges are allowed to be flipped. (It is not clear what it means to
flip an edge of the outerface since for $n>3$, the outerface is not a
triangle.)\ A flip in an outerplanar graph corresponds to a certain rotation in the dual tree; see \citep{CCW-IJCM05, Pallo-IPL00, Pallo-IPL87, STT-JAMS88}. This section focuses on simulatenous flips in maximal outerplane graphs, which have not previously been studied.

\begin{lemma} 
\lemlabel{OuterplaneFlip}
Every internal edge of a maximal outerplane graph is flippable.  
\end{lemma}

\begin{proof} 
Suppose that an internal edge $vw$ is not filppable. Then $vw$ is blocked by some edge $pq$. Thus $\{v,w,p,q\}$ induce $K_4$. This is a contradiction since no outerplane graph contains $K_4$.
\end{proof}

%\begin{proof}[Old Proof] Suppose that an internal edge $vw$ is blocked by some edge $pq$. Without loss of generality, the vertices appear in the order $(v,p,w,q)$ on the outerface. Thus the boundary of the outerface along with $vw$ and $pq$ form a subdivision of $K_4$. This is a contradiction, since no outerplane graph contains a subdivision of $K_4$. Thus $vw$ is not blocked and is flippable. \end{proof}

\begin{lemma}
\lemlabel{OuterplaneFlippable}
A set $S$ of internal edges in a maximal outerplane graph $G$ is flippable if and only if the corresponding dual edges $S^*$ form a matching in $G^*$.
\end{lemma}

\begin{proof}
For $S$ to be flippable it is necessary that there are no two consecutive edges in $S$. This is equivalent to the condition that $S^*$ is a matching of $G^*$. 
By \lemref{OuterplaneFlip}, every edge in $S$ is flippable. As in \lemref{FlippableCharacterisation}, the only obstruction to $S$ being flippable is a bad pair of edges, which cannot occur since a bad pair of edges contains a subdivision of $K_4$.
\end{proof}

\begin{theorem}
\thmlabel{OuterplaneBigFlip}
Every $n$-vertex maximal outerplane graph $G$ has a flippable set of at least $\third(n-3)$ edges. Moreover, for infinitely many $n$, there is an $n$-vertex maximal outerplane graph in which every flippable set has at most $\third(n-3)$ edges. 
\end{theorem}

\begin{proof}
First we prove the lower bound. Since $G^*$ is a tree with maximum degree at most three, $G^*$ has a proper edge $3$-colouring (by an easy inductive argument). Now $G^*$ has $n-3$ edges. Thus the largest colour class is a matching of at least $\third(n-3)$ edges, which by \lemref{OuterplaneFlippable}, corresponds to a flippable set of at least $\third(n-3)$ edges in $G$.

Now we prove the upper bound. By \lemref{OuterplaneFlippable} it suffices to construct trees $T$ with maximum degree three, in which the maximum cardinality of a matching equals one third the number of edges. We can then take the  maximal outerplane graph $G$ for which $G^*=T$. Let $T$ be a tree rooted at a vertex $r$ such that every non-leaf vertex has degree three, and the distance between every leaf vertex and $r$ is odd. Obviously there are infinitely many such trees. Let $K$ be the set of vertices at even distance from $r$. Then $K$ is a \emph{vertex cover} of $T$ (that is, every edge of $T$ is  incident to a vertex in $K$). Since no edge of $T$ has both its endpoints in $K$, and every vertex in $K$ has degree three, $|K|$ equals one third the number of edges. Since $T$ has maximum degree three, $K$ is a minimum vertex cover. \citet{Konig36} proved that the maximum cardinality of a matching in a bipartite graph equals the minimum cardinality of a vertex cover. Thus the maximum cardinality of a matching equals one third the number of edges, as desired.
\end{proof}

%%%%%%%%%%%%%%%%%%%%%%%%%%%%%%%%%%%%%%%%%%%%%%%%%%%%%%%%%%%%%%%%%%%%%%%%%%%%%

The following is the main result of this section. In the remainder of this paper all logarithms have base $2$, and $c_1$ is the constant $2/\log\tfrac{6}{5}$ ($\approx 7.6$). 

\begin{theorem}
\thmlabel{TwoOuterplane}
Let $G_1$ and $G_2$ be (unlabelled) maximal outerplane graphs on $n$ vertices. There is a sequence of at most $4c_1\log n$ simultaneous flips to transform $G_1$ into $G_2$.
\end{theorem}

\thmref{TwoOuterplane} is implied by the following lemma. 

\begin{lemma}
\lemlabel{MakeOuterplaneDominant}
For every maximal outerplane graph $G$ on $n$ vertices, there is a sequence of at most $2c_1\log n$ simultaneous flips to transform $G$ into a maximal outerplane graph that has a dominant vertex.
\end{lemma}

\begin{proof}[Proof of \thmref{TwoOuterplane} assuming
\lemref{MakeOuterplaneDominant}] Observe that two $n$-vertex maximal
outerplane graphs each with a dominant vertex are isomorphic. Let $D_n$ denote
the $n$-vertex maximal outerplane graph with a dominant vertex. To
transform $G_1$ into $G_2$, first transform $G_1$ into $D_n$, and then apply
the flips to transform $G_2$ into $D_n$ in reverse order.  
\end{proof}

The proof of  \lemref{MakeOuterplaneDominant} proceeds in two parts. In \lemref{OuterplaneDiameter} we reduce the diameter of the dual tree to $c_1\log n$ using at most $c_1\log n$ simultaneous flips. Then in \lemref{OuterplaneDominant} a dominant vertex is introduced using at most a further $c_1\log n$ simultaneous flips.

%%%%%%%%%%%%%%%%%%%%%%%%%%%%%%%%%%%%%%%%%%%%%%%%%%%%%%%%%%%%%%%%%%%%%%%%%%%%%%%

\begin{lemma}  
\lemlabel{OuterplaneDiameter}
Let $G$ be a maximal outerplane graph on $n$ vertices. Then $G$ can be transformed by a sequence of at most $c_1\log n$ simultaneous flips into a maximal outerplane graph $X$ such that the diameter of the dual tree $X^*$ is at most $c_1\log n$.
\end{lemma}

\begin{proof} We proceed by induction on $n$. The result holds trivially for $n=3$. Assume the lemma holds for graphs with less than $n$ vertices, and let $G$ be a maximal outerplane graph on $n$ vertices. By a theorem of \citet{BDW-CDM06}, $G$ has an independent set $I$ of at least $\tfrac{n}{6}$ vertices, and $\deg_G(v)\leq4$ for every vertex $v\in I$. Obviously $\deg_G(v)\geq2$. For $d\in\{2,3,4\}$, let $I_d:=\{v\in I:\deg_G(v)=d\}$. 

For every vertex $v\in I_3\cup I_4$, add one internal edge incident to $v$ to a set $S$. Since $I$ is independent,  $|S|=|I_3|+|I_4|$. Suppose on the contrary that there are two consecutive edges $xu,xv\in S$. Then $x\not\in I_3\cup I_4$, which implies that $u,v\in I_3\cup I_4$. Since every internal face of $G$ is a triangle, $uv$ is an edge of $G$, which contradicts the independence of $I$. Thus no two edges in $S$ are consecutive. By \lemref{OuterplaneFlippable}, $S$ is flippable in $G$. Let $G':=\F{G}{S}$. Every vertex $v\in I_2\cup I_3$ has $\deg_{G'}(v)=2$, and every vertex $v\in I_4$ has $\deg_{G'}(v)=3$. 

Since $I_4$ is an independent set of $G$, and any edge in $G'$ that is incident to a vertex in $I_4$ is also in $G$, $I_4$ is an independent set of $G'$. Let $S'$ be the set of internal edges of $G'$ incident to a vertex in $I_4$. Thus $|S'|=|I_4|$, and by the same argument used for $S$, no two edges in $S'$ are consecutive in $G'$. By \lemref{OuterplaneFlippable}, $S'$ is flippable in $G'$. Let $G'':=\F{G'}{S'}$. Every vertex $v\in I$ has $\deg_{G''}(v)=2$. 

Thus $G$ can be transformed by two simultaneous flips into a maximal outerplane graph $G''$ containing at least $\tfrac{n}{6}$ vertices of degree two. Let $G'''$ be the  maximal outerplane graph obtained from $G''$ by deleting the vertices of degree two. Then $G'''$ has at most $\tfrac{5}{6}n$ vertices. By induction, $G'''$ can be transformed by a sequence of at most $c_1\log\tfrac{5}{6}n$ simultaneous flips into a maximal outerplane graph $X$ such that the diameter of $X^*$ is at most $c_1\log\tfrac{5}{6}n$. 

Consider a vertex $v\in I$. Since $\deg_{G''}(v)=2$, there is one internal face incident to $v$ in $G''$, which corresponds to a leaf in $G''^*$. Thus the dual tree $X^*$ is obtained by adding leaves to the dual tree $G''^*$. Hence the diameter of $X^*$ is at most the diameter of $G''^*$ plus $2$, which is $2+c_1\log\tfrac{5}{6}n=c_1\log n$. We have used two simultaneous flips, $S$ and $S'$, to transform $G$ into $G''$, and then $c_1\log\tfrac{5}{6}n$ simultaneous flips to transform $G''$ into $X$. The total number of flips is $2+c_1\log\tfrac{5}{6}n=c_1\log n$. \end{proof}

%%%%%%%%%%%%%%%%%%%%%%%%%%%%%%%%%%%%%%%%%%%%%%%%%%%%%%%%%%%%%%%%%%%%%%%%%%%

\begin{lemma} 
\lemlabel{OuterplaneDominant}
Let $G$ be a maximal outerplane graph on $n$ vertices. Suppose that $G^*$ has diameter $k$. Let $v$ be a fixed vertex of $G$. Then $G$ can be transformed by at most $k$ simultaneous flips into a maximal outerplane graph $X$ in which $v$ is dominant. 
\end{lemma}

\begin{proof} 
As illustrated in \figref{MakeDominantOuterplane}, let $P$ be the set of internal faces incident with $v$ in $G$. In the dual tree $G^*$, the corresponding vertices of $P$ form a path $P^*$. Define the \emph{distance} of each vertex $x$ in $G^*$ as the minimum number of edges in a path from $x$ to a vertex in $P^*$. Since the diameter of $G^*$ is $k$, every vertex in $G^*$ has distance at most $k$. No two vertices in $G^*$ both with distance one are adjacent, as otherwise $G^*$ would contain a cycle. Each vertex of $P^*$ is adjacent to at most one vertex at distance one, since $G^*$ has maximum degree at most three, and the endpoints of $P^*$ correspond to faces with an edge on the outerface of $G$. Let $S^*$ be the set of edges of $G^*$ incident to $P^*$ but not in $P^*$. Then $S^*$ is a matching between the vertices at distance one and the vertices of $P^*$, such that all vertices at distance one are matched. Let $S$ be the set of edges of $G$ corresponding to $S^*$ under duality. 
Note that $S$ is the set of internal edges that are seen by $v$. By \lemref{OuterplaneFlippable}, $S$ is a flippable set of edges of $G$. Let $G':=\F{G}{S}$. In $G'$, the distance of each vertex not adjacent to $P^*$ is reduced by one. Thus, by induction, at most $k$ simultaneous flips are required to reduce the distance of every vertex to zero, in which case $v$ is adjacent to every other vertex.
%, and every edge incident to $v$ is in $G$. WHY DID I HAVE THIS?
\end{proof}

\Figure{MakeDominantOuterplane}{\includegraphics{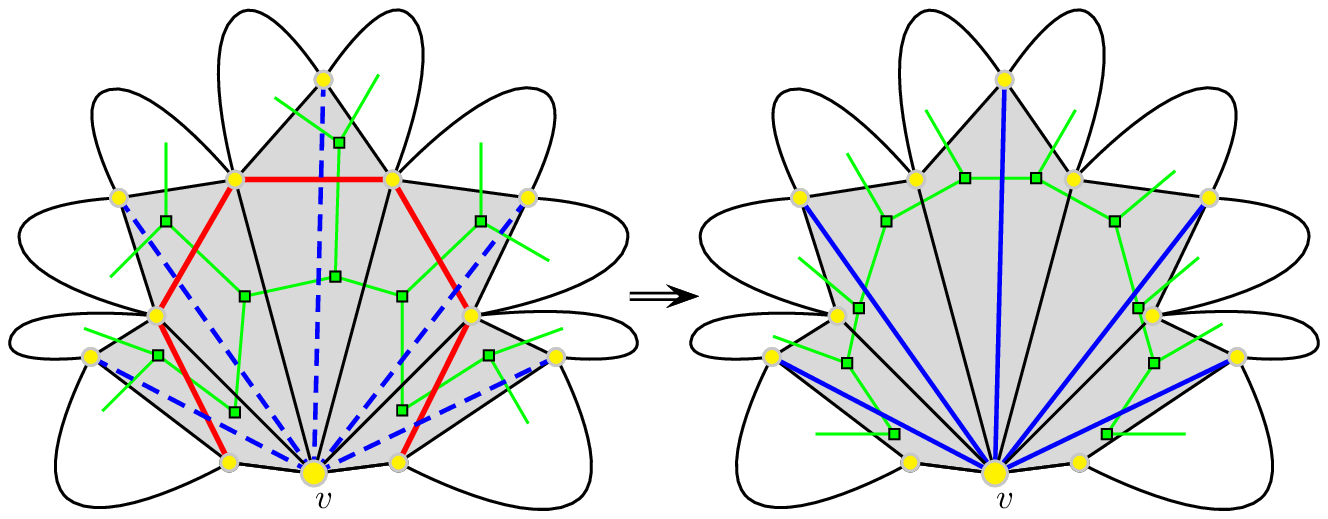}}{Making $v$ a dominant vertex in \lemref{OuterplaneDominant}; the vertices of the dual tree are drawn as squares.}

Clearly \lemref{MakeOuterplaneDominant} is implied by \twolemref{OuterplaneDiameter}{OuterplaneDominant} (with $k=c_1\log n$).

%%%%%%%%%%%%%%%%%%%%%%%%%%%%%%%%%%%%%%%%%%%%%%%%%%%%%%%%%%%%%%%%%%%%%%%%%%%%%%
\section{Simultaneous Flips Between Given Triangulations}
\seclabel{TwoTriangulations}
%%%%%%%%%%%%%%%%%%%%%%%%%%%%%%%%%%%%%%%%%%%%%%%%%%%%%%%%%%%%%%%%%%%%%%%%%%%%%%

In this section we prove the following theorem, which is an analogue of 
\thmref{TwoOuterplane} for outerplane graphs. Throughout, $c_1$ is the constant $2/\log\tfrac{6}{5}$ ($\approx 7.6$) from \secref{Outerplane}, and $c_2$ is the constant $2/\log\tfrac{54}{53}$ ($\approx74.2$).

\begin{theorem}
\thmlabel{TwoTriangulations}
Let $G_1$ and $G_2$ be (unlabelled) triangulations on $n$ vertices. There is a sequence of at most $2+4(c_1+c_2)\log n$ simultaneous flips to transform $G_1$ into $G_2$.
\end{theorem}

\thmref{TwoTriangulations} is implied by the following lemma using the approach of Wagner described in \secref{Introduction}.

\begin{lemma}
\lemlabel{MakeDoubleDominant}
For every $n$-vertex triangulation $G$, there is a sequence of at most $1+2(c_1+c_2)\log n$ simultaneous flips to transform $G$ into the standard triangulation $\Delta_n$.
\end{lemma}

%\begin{proof}[Proof of \thmref{TwoTriangulations} assuming \lemref{MakeDoubleDominant}] To transform $G_1$ into $G_2$, first transform $G_1$ into $D_n$, and then apply the flips to transform $G_2$ into $D_n$ in reverse order.\end{proof}

To prove \lemref{MakeDoubleDominant} we first apply \thmref{Hamiltonian} to obtain a Hamiltonian triangulation with one simultaneous flip. Thus it suffices to prove that a Hamiltonian triangulation can be transformed into $\Delta_n$. A Hamiltonian cycle $H$ of a triangulation $G$ naturally divides $G$ into two maximal outerplane subgraphs: an `inner' subgraph consisting of $H$ and the edges inside $H$, and an `outer' subgraph consisting of $H$ and the edges outside of $H$. (Note that \citet{MNO-GC03} used a similar approach for individual flips.)\ At this point, it is tempting to apply \lemref{MakeOuterplaneDominant} twice, once on the inner subgraph to obtain one dominant vertex, and then on the outer subgraph to obtain a second dominant vertex, thus reaching the standard triangulation. However, \lemref{MakeOuterplaneDominant} cannot be applied directly since we need to take into consideration the interaction between these two outerplane subgraphs. The main problem is that an internal edge in the inner subgraph may be blocked by an edge in the outer subgraph. The bulk of this section is dedicated to solving this impasse. 

First some definitions. A \emph{chord} of a cycle $C$ in a triangulation $G$ is an edge of $G$ that is not in $C$ and whose endpoints are both in $C$. A chord $e$ of $C$ is classified as \emph{internal} or \emph{external} depending on whether $e$ is contained in the interior or exterior of $C$ (with respect to the outerface of $G$). For the inductive step in \lemref{InternalDiameter} below to work we need to consider a more general type of cycle than a Hamiltonian cycle. A cycle $C$ of a triangulation $G$ is \emph{empty} if the interior of $C$ contains no vertices of $G$. Obviously a Hamiltonian cycle is always empty. For an empty cycle $C$ of a triangulation $G$, let $G\{C\}$ denote the subgraph of $G$ whose vertices are the vertices of $C$, and whose edges are the edges of $C$ along with the internal chords of $C$. Then $G\{C\}$ is a maximal outerplane graph, and the boundary of the outerface of $G\{C\}$ is $C$. 

\begin{lemma}
\lemlabel{InternalIsFlippable}
Let $C$ be an empty cycle of a triangulation $G\ne K_4$. Let $vw$ be an internal chord of $C$ that is blocked by some edge $pq$. Then $pq$ is an external chord of $C$ that is flippable in $G$.
\end{lemma}

\begin{proof}
By \lemref{BlockingEdgeFlippable}, $pq$ is a flippable edge of $G$. Since $C$ is empty, $p$ and $q$ are vertices of $C$. Now $pq$ is not internal, as otherwise $\{p,q,v,w\}$ would induce $K_4$ in the outerplane graph $G\{C\}$. Thus $pq$ is external. 
\end{proof}

%Without loss of generality $v$ is in the interior of the region bounded by $C$ and the edge $pq$. Let $x$ and $y$ be the vertices seen by $pq$. Without loss of generality, $x=v$ or $x$ is in the interior of the triangle $(p,v,q)$. First suppose that $x\ne v$. Then $(p,v,q)$ is a separating triangle, since $w$ is in its exterior. Thus $pq$ is flippable by \lemref{EdgeInSepTriangle}. Now suppose that $x=v$. We have that $y\ne w$ as otherwise $G=K_4$. Thus $v$ is in the interior of the triangle $(p,w,q)$ and $y$ is in the exterior. Hence $(p,w,q)$ is a separating triangle, and $pq$ is flippable by \lemref{EdgeInSepTriangle}. 

\begin{lemma}
\lemlabel{FlipInternalMatching}
Let $C$ be an empty cycle of a triangulation $G$. Let $S$ be a set of internal chords of $C$, no two of which are consecutive. Then there is a flippable set $T$ of edges in $G$ such that:\\
\aaa\ $T\cap C=\emptyset$,\\
\bbb\ $|S\cap T|\geq\third|S|$, and\\
\ccc\ every edge in $T\setminus S$ is an external chord of $C$ and  $|T\setminus S|\leq|S\cap T|$.
\end{lemma}

\begin{proof}
Let $S'$ be the set of edges in $S$ that are individually flippable in $G$. Let $S'':=S\setminus S'$. By \lemref{InternalIsFlippable}, there is an external chord that blocks each edge $e\in S''$. Distinct edges $e_1,e_2\in S''$ are blocked by distinct external chords, as otherwise $e_1$ and $e_2$ would be a bad pair, and the outerplane graph $G\{C\}$ would contain a subdivision of $K_4$. Let $B$ be this set of blocking external chords. Thus $|B|=|S''|$. By \lemref{TaitColouring}, $B$ can be $3$-coloured such that no two monochromatic edges in $B$ are consecutive in $G$. (Note that since $B$ forms an outerplane subgraph of $G$, this $3$-colouring can be computed in \Oh{n} time without using the $4$-colour theorem.)\ Let $P$ be the largest set of monochromatic edges in $B$. Then $|P|\geq\third|B|$. Let $Q$ be the set of edges in $S''$ that are blocked by edges in $P$. Then $|Q|=|P|$. Let $T:=S'\cup P\cup Q$. Observe that $T\cap C=\emptyset$. This proves (a).

To prove that $T$ is flippable in $G$, we verify each of the conditions of \lemref{FlippableCharacterisation}. $T$ consists of internal chords $S'\cup Q$, and external chords $P$. Since $S'\cup Q\subseteq S$, no two internal chords in $T$ are consecutive. By the construction of $P$, no two external chords in $T$ are consecutive. Since the internal chords and external chords are separated by $C$, no two edges in $T$ are consecutive. Thus condition (1) of \lemref{FlippableCharacterisation} is satisfied. 

As in \lemref{OuterplaneFlippable}, there is no bad pair among the internal chords as otherwise $G\{C\}$ would contain a subdivision of $K_4$. Similarly there is no bad pair among the external chords. Suppose there is a bad pair of edges in $T$, one an internal chord $xy$ and the other an external chord $vw$. Then both $vw$ and $xy$ are seen by some pair of vertices $p$ and $q$. Since $vw\in P\subseteq B$, $vw$ blocks some internal chord $ab\in S''$. By \lemref{BlockingEdgeInBadPair}, $ab$ and $xy$ are consecutive, which is a contradiction since both edges are in $S$. Thus there is no bad pair in $T$, and condition (2) of \lemref{FlippableCharacterisation} is satisfied.

Each edge in $P$ blocks some other edge, and is thus individually flippable by \lemref{BlockingEdgeFlippable}. By definition, all the edges in $S'$ are individually flippable in $G$. While each edge in $Q$ is not individually flippable, the corresponding blocking edge is in $P\subseteq T$. Thus condition (3) of \lemref{FlippableCharacterisation} is satisfied. Therefore $T$ is flippable in $G$.

Now $T\cap S=S'\cup Q$. Since $S'\cap Q=\emptyset$, we have $|T\cap S|=|S'|+|Q|=|S'|+|P|\geq|S'|+\third|B|\geq\third|S'|+\third|S''|=\third|S|$. This proves (b). Now $T\setminus S=P$, all of whose elements are external chords. Since $|S\cap T|=|S'|+|P|$, we have $|P|\leq|S\cap T|$. Since $T\setminus S=P$, we have  $|T\setminus S|\leq|S\cap T|$. This proves (c).
\end{proof}

%%%%%%%%%%%%%%%%%%%%%%%%%%%%%%%%%%%%%%%%%%%%%%%%%%%%%%%%%%%%%%%%%%%%%%%%%%%%%

The following result extends \lemref{OuterplaneDiameter} for outerplane graphs to the case of triangulations.

%The \emph{internal diameter} of an empty cycle $C$ in a triangulation $G$ is the diameter of $G\{C\}^*$. The first step in our algorithm to transform a Hamiltonian triangulation into the standard triangulation is to reduce the internal diameter of the Hamiltonian cycle. 

\begin{lemma}  
\lemlabel{InternalDiameter}
Let $H$ be a Hamiltonian cycle of a triangulation $G$ with $n$ vertices. Then $G$ can be transformed by a sequence of at most $c_2\log n$ simultaneous flips into a triangulation $X$ in which $H$ is a Hamiltonian cycle and the diameter of $X\{H\}^*$ is at most $c_2\log n$.
\end{lemma}

\begin{proof} We proceed by induction on $n$ with the following stronger hypothesis:

``Let $G$ be a triangulation, and let $C$ be an empty cycle of $G$ with $n$ vertices. ($G$ may have more than $n$ vertices.)\ Then $G$ can be transformed by a sequence of at most $c_2\log n$ simultaneous flips into a triangulation $X$ in which $C$ is an empty cycle and the diameter of $X\{C\}$ is at most $c_2\log n$. Moreover, every edge of $G$ that is incident to a vertex not in $C$ remains in $X$.'' 

The lemma immediately follows since any Hamiltonian cycle is empty. The hypothesis holds trivially for $n=3$. Assume the hypothesis holds for all triangulations with less than $n$ vertices. Let $G$ be a triangulation, and let $C$ be an empty cycle of $G$ with $n$ vertices. 

By a theorem of \citet{BDW-CDM06}, the outerplane graph $G\{C\}$ has an independent set $I$ of at least $\tfrac{n}{6}$ vertices, and $\deg_{G\{C\}}(v)\leq4$ for every vertex $v\in I$. Obviously $\deg_{G\{C\}}(v)\geq2$. For $d\in\{2,3,4\}$, let $I_d:=\{v\in I:\deg_{G\{C\}}(v)=d\}$. 

For every vertex $v\in I_3\cup I_4$, add one internal chord of $C$ that is incident to $v$ to a set $S$. Since $I$ is independent, $|S|=|I_3|+|I_4|$.  Suppose on the contrary that there are two consecutive edges $xu$ and $xv$ in $S$. Then $x\not\in I_3\cup I_4$, which implies that $u,v\in I_3\cup I_4$. Since every face of $G$ is a triangle, $uv$ is an edge, which contradicts the independence of $I$. Thus no two edges in $S$ are consecutive. By  \lemref{FlipInternalMatching}, there is a flippable set of edges $T$ in $G$, such that  $T\cap C=\emptyset$ and $|S\cap T|\geq\third|S|=\third(|I_3|+|I_4|)$. Moreover, every edge in $T\setminus S$ is an external chord of $C$ in $G$. For $d\in\{3,4\}$, let $I_d'$ be the set of vertices in $I_d$ incident to an edge in $S\cap T$. Thus $|I_3'|+|I_4'|\geq\third(|I_3|+|I_4|)$.

Let $G':=\F{G}{T}$. Since $T\cap C=\emptyset$, $C$ is an empty cycle of $G'$. Every vertex $v\in I_2\cup I'_3$ has $\deg_{G'\{C\}}(v)=2$. Every vertex $v\in I'_4$ has $\deg_{G'\{C\}}(v)=3$. 

An edge in $G'\{C\}$ that is incident to a vertex in $I_4'$ is also in $G\{C\}$. Since $I_4'$ is an independent set of $G\{C\}$, it is also an independent set of $G'\{C\}$. Let $S'$ be the set of internal chords of $C$ in $G'$ that are incident to a vertex in $I_4'$. Thus $|S'|=|I_4'|$, and by the same argument used for $S$, no two edges in $S'$ are consecutive in $G'$. By  \lemref{FlipInternalMatching}, there is a flippable set of edges $T'$ in $G'$, such that  $T'\cap C=\emptyset$ and $|S'\cap T'|\geq\third|S'|=\third|I_4'|$. Moreover, every edge in $T'\setminus S'$ is an external chord of $C$ in $G'$.  Let $I_4''$ be the set of vertices in $I_4'$ incident to an edge in $S'\cap T'$. Thus $|I_4''|\geq\third|I_4'|$.

Let $G'':=\F{G'}{T'}$. Since $T'\cap C=\emptyset$, $C$ is an empty cycle of $G''$. Every vertex $v\in I_2\cup I_3'\cup I''_4$ has $\deg_{G''\{C\}}(v)=2$. Now 
$|I_2\cup I_3'\cup I''_4|
\geq|I_2|+|I_3'|+\third|I_4'|
\geq|I_2|+\third(|I_3'|+|I_4'|)
\geq|I_2|+\ninth(|I_3|+|I_4|)
\geq\ninth(|I_2|+|I_3|+|I_4|)
=\ninth|I|
\geq\tfrac{n}{54}$.
  
In summary, $G$ can be transformed by two simultaneous flips into a triangulation $G''$ in which $C$ is an empty cycle, and $G''\{C\}$ has an independent set $L$ ($=I_2\cup I_3'\cup I''_4$) such that $|L|\geq\tfrac{n}{54}$ and $\deg_{G''\{C\}}(v)=2$ for every vertex $v\in L$. Consider a vertex $v\in L$. Say $(u,v,w)$ is the $2$-edge path in $C$. Since $L$ is independent, $u\not\in L$ and $w\not\in L$. Since $\deg_{G''\{C\}}(v)=2$, $uw$ is an internal chord of $C$ in $G''$. Let $D$ be the cycle of $G$ obtained by replacing the the path $(u,v,w)$ in $C$ by the edge $uw$ (for all $v\in L$). Thus $D$ is an empty cycle of $G''$, and $|D|=n-|L|\leq\tfrac{53}{54}n$. By induction applied to $D$ and $G''$, $G''$ can be transformed by a sequence of at most $c_2\log\tfrac{53}{54}n$ simultaneous flips into a triangulation $X$ in which $D$ is an empty cycle and the diameter of $X\{D\}^*$ is at most $c_2\log\tfrac{53}{54}n$. Moreover, every edge of $G''$ that is incident to a vertex not in $D$ remains in $X$. 

Consider a vertex $v\in L$. Say $(u,v,w)$ is the $2$-edge path in $C$. Since $v$ is not in $D$, the edges $uv$ and $vw$ of $G$ are in $X$. Thus $C$ is an empty cycle of $X$. Since $uw$ is an edge of $D$, $uvw$ is a face of $X$. The vertex in the dual tree $X\{C\}^*$ that corresponds to $uvw$ is a leaf in $X\{C\}^*$. Thus $X\{C\}^*$ is obtained by adding leaves to the dual tree $X\{D\}^*$. Hence the diameter of $X\{C\}^*$ is at most the diameter of $X\{D\}^*$ plus $2$, which is at most $2+c_2\log\tfrac{53}{54}n=c_2\log n$. We have used two simultaneous flips, $T$ and $T'$, to transform $G$ into $G''$, and then at most  $c_2\log\tfrac{53}{54}n$ simultaneous flips to transform $G''$ into $X$. The total number of flips is at most $2+c_2\log\tfrac{53}{54}n=c_2\log n$. Since every edge in $T$ is a chord of $C$ in $G$, and every edge in $T'$ is a chord of $C$ in $G'$, every edge of $G$ that is incident to a vertex not in $C$ remains in $X$. \end{proof}

%%%%%%%%%%%%%%%%%%%%%%%%%%%%%%%%%%%%%%%%%%%%%%%%%%%%%%%%%%%%%%%%%%%%%%%%%%%

The following result is analogous to \lemref{OuterplaneDominant} for outerplane graphs. The key difference is that the choice of vertex $v$ is no longer arbitrary.

\begin{lemma}
\lemlabel{Star}  
Let $H$ be a Hamiltonian cycle of a triangulation $G$. Suppose that $G\{H\}^*$ has diameter $k$. Let $v$ be a vertex of $G$ not incident to any external chords of $H$ in $G$. Then $G$ can be transformed  by at most $k$ simultaneous flips
into a triangulation $X$ in which $H$ is a Hamiltonian cycle of $X$ and $v$ is dominant. Moreover, every edge incident to $v$ is in $X\{H\}$.
\end{lemma}

\begin{proof} First note that there is such a vertex $v$ since the subgraph of $G$ consisting of $H$ and the external chords of $H$ is maximal outerplane, and thus has a vertex of degree two. Let $P$ be the set of internal faces incident with $v$ in $G$. In the dual tree $G\{H\}^*$, the corresponding vertices of $P$ form a path $P^*$. Define the \emph{distance} of each vertex $x$ in $G\{H\}^*$ as the minimum number of edges in a path from $x$ to a vertex in $P^*$. Since the diameter of $G\{H\}^*$ is $k$, every vertex in $G^*$ has distance at most $k$. No two vertices in $G\{H\}^*$ both with distance one are adjacent, as otherwise $G\{H\}^*$ would contain a cycle. Each vertex of $P^*$ is adjacent to at most one vertex at distance one, since $G\{H\}^*$ has maximum degree at most three, and the endpoints of $P^*$ correspond to faces with an edge on the outerface of $G\{H\}$. Let $S^*$ be the set of edges of $G\{H\}^*$ incident to $P^*$ but not in $P^*$. Then $S^*$ is a matching between the vertices at distance one and the vertices of $P^*$, such that all vertices at distance one are matched. Let $S$ be the set of edges of $G\{H\}$ corresponding to $S^*$ under duality. Consider an edge $xy\in S$. Then $xy$ is seen by $v$ and some other vertex $w$. If $xy$ is not flippable, then by \lemref{InternalIsFlippable}, $vw$ is an external chord of $H$ in $G$. Thus $xy$ is flippable, since by construction, $v$ is not incident to any external chords of $H$ in $G$. Hence $S$ is a set of individually flippable edges. No two edges in $S$ are consecutive, since every internal face of $G\{H\}$ is a triangle. No two edges in $G\{H\}$ form a bad pair since $G\{H\}$ is outerplane. By \lemref{FlippableCharacterisation}, $S$ is flippable in $G$. Let $G':=\F{G}{S}$. Observe that $S\cap H=\emptyset$. Thus $H$ is a Hamiltonian cycle of $G'$. In $G'\{H\}$, the distance of each vertex not adjacent to $P^*$ is reduced by one. Thus, by induction, at most $k$ simultaneous flips are required to reduce the distance of every vertex to zero, in which case $v$ is adjacent to every other vertex of $G$, and every edge incident to $v$ is in $G\{H\}$. \end{proof}

\twolemref{InternalDiameter}{Star} imply:

\begin{lemma} 
\lemlabel{PartOne}
Let $H$ be a Hamiltonian cycle of a triangulation $G$. Then $G$ can be transformed  by at most $2c_2\log n$ simultaneous flips into a triangulation $X$ in which $H$ is a Hamiltonian cycle of $X$, and there is a vertex $v$ adjacent to every other vertex, and every edge incident to $v$ is in $X\{H\}$.\qed
\end{lemma}

We are now half way to transforming a given triangulation into the standard triangulation. The second half is somewhat easier.

\begin{lemma}
\lemlabel{PartTwo}
Let $G$ be an $n$-vertex triangulation with a dominant vertex $v$. Then there is a sequence of at most $2c_1\log(n-1)$ simultaneous flips to transform $G$ into the standard triangulation on $n$ vertices.
\end{lemma}

\begin{proof}
Observe that $G\setminus v$ is a maximal outerplane graph, in which the vertices are ordered on the outerface according to the cyclic order of the neighbours of $v$. Let $C$ be the cycle bounding the outerface of $G\setminus v$. By \lemref{MakeOuterplaneDominant} there is a sequence of at most $2c_1\log(n-1)$ simultaneous flips to transform $G\setminus v$ into a maximal outerplane graph with a dominant vertex. Each of these flips is valid in $G$ since $C$ has no internal chords (cf.~\lemref{FlipInternalMatching}). We obtain the standard triangulation.
\end{proof}

Observe that \twolemref{PartOne}{PartTwo} together prove \lemref{MakeDoubleDominant}, which in turn proves \thmref{TwoTriangulations}. Although the \Oh{\log n} simultaneous flips in \thmref{TwoTriangulations} may each involve a linear number of edges, the total number of flipped edges is linear. 

\begin{theorem}
\thmlabel{TwoTriangulationsTotalFlips}
Let $G_1$ and $G_2$ be triangulations on $n$ vertices. There is a sequence of \Oh{\log n} simultaneous flips to transform $G_1$ into $G_2$, and \Oh{n} edges are flipped in total.
\end{theorem}

\begin{proof} 
It suffices to prove that there are \Oh{n} flips in \twolemref{InternalDiameter}{Star}, since  at most $n$ edges are flipped to make the graph Hamiltonian, and there are constant times as many  flips in \thmref{TwoTriangulations} as there are in \twolemref{Star}{InternalDiameter}. 
In \lemref{Star}, each flipped edge becomes incident to $v$, and then remains incident to $v$. Thus the number of flipped edges is at most $n-1$. In \lemref{InternalDiameter}, \Oh{n} edges are flipped to obtain a triangulation on at most $\frac{53}{54}n$ vertices. Therefore, the number of flipped edges $F(n)$ satisfies the recurrence  $F(n)=F(\frac{53}{54}n)+\Oh{n}$, which solves to \Oh{n}. 
\end{proof}

%%%%%%%%%%%%%%%%%%%%%%%%%%%%%%%%%%%%%%%%%%%%%%%%%%%%%%%%%%%%%%%%%%%%%%%%%%%%%%
\section{Large Simultaneous Flips}
\seclabel{BigFlip}
%%%%%%%%%%%%%%%%%%%%%%%%%%%%%%%%%%%%%%%%%%%%%%%%%%%%%%%%%%%%%%%%%%%%%%%%%%%%%%

In this section we prove bounds on the size of a maximum simultaneous flip in a
triangulation. Let $\msf{G}$  denote the maximum cardinality of a flippable set of edges in a triangulation $G$. In related work, \citet{GUW-GC01} proved that every triangulation has at least $n-2$ (individually) flippable edges, and every triangulation with minimum degree four has at least $2n+3$ (individually)
flippable edges. \citet{GHNP-IJCGA03} proved that every geometric triangulation has a set of at least $\sixth(n-4)$ simultaneously flippable edges.
%
%\subsection{Lower Bound}
%%%%%%%%%%%%%%%%%%%%%%%%%%%%%%%%%%%%%%%%%%%%%%%%%%%%%%%%%%%%%%%%%%%
%
The following is the main result of this section.

\begin{theorem}
\thmlabel{BigFlip}
For every triangulation $G$ with $n\geq4$ vertices, $\msf{G}\geq\third(n-2)$.
\end{theorem}

Assume there is a \emph{counterexample} to \thmref{BigFlip}; that is,  a triangulation $G$ with $n\geq4$ vertices and $\msf{G}<\third(n-2)$. A counterexample with the minimum number of vertices is a \emph{minimum counterexample}. 

\begin{lemma}
\lemlabel{SmallTriangBigFlip}
A counterexample has $n\geq 7$ vertices.
\end{lemma}

\begin{proof}
If $n=4$ then $G=K_4$, which has a flippable set of $2>\third(4-2)$ edges, as illustrated in \figref{FourFiveVertex}(a). If $n=5$ then $G=K_5\setminus e$, which has a flippable set of $2>\third(5-2)$ edges. 

\Figure{FourFiveVertex}{\includegraphics{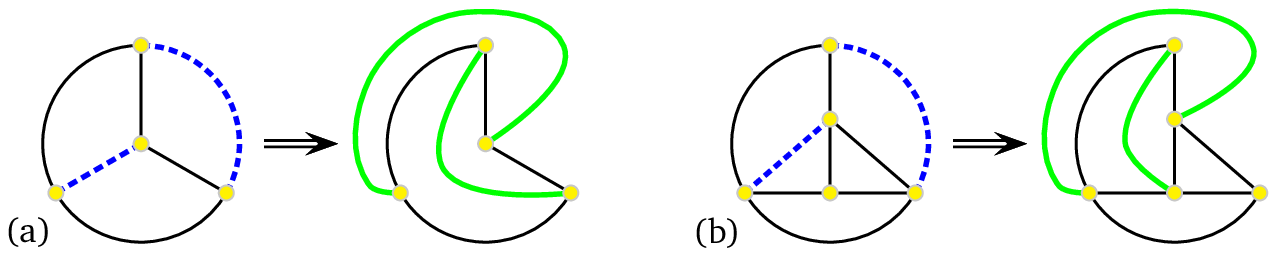}}{2-edge flip in (a) $K_4$ and (b) $K_5\setminus e$.}

If $n=6$ then $G$ is the octahedron illustrated in \figref{SixVertex}(a), or $G$ is the triangulation illustrated in \figref{SixVertex}(b). In both cases there is a flippable set of $3>\third(6-2)$ edges. 
\end{proof}

\Figure{SixVertex}{\includegraphics{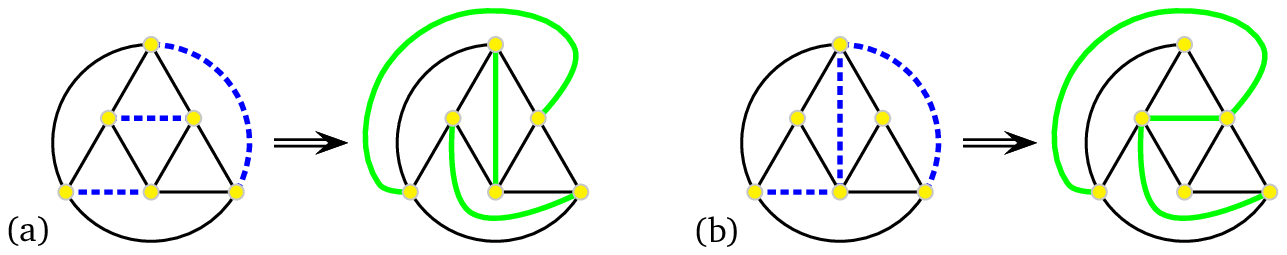}}{$3$-edge flip in (a) the octahedron and (b) the other $6$-vertex triangulation.}

%%%%%%%%%%%%%%%%%%%%%%%%%%%%%%%%%%%%%%%%%%%%%%%%%%%%%%%%%%%%%%%%%%%%%%%%%%%%%%

\begin{lemma}
\lemlabel{ThreeFour}
A minimum counterexample has no edge $vw$ with $\deg(v)=3$ and $\deg(w)=4$.
\end{lemma}

\begin{proof} 
Let $G$ be a minimum counterexample with $n$ vertices. Suppose that $G$ has an edge $vw$ with $\deg(v)=3$ and $\deg(w)=4$. Then the neighbours of $v$ and $w$ form a triangle $(x,y,z)$ with $v$ adjacent to $x$ and $y$, and $w$ adjacent to $x$, $y$ and $z$. Let $G':=(G\setminus v)\setminus w$. Then $G'$ is a triangulation with $n-2$ vertices in which $(x,y,z)$ is a face. Since $G$ is minimum, $G'$ is not a counterexample. Thus $G'$ has a flippable set $S'$ of at least $\third(n-4)$ edges. At most one of $\{xy,xz,yz\}$ is in $S'$. If $xz\in S'$, then let $S:=S'\cup\{yw\}$. Otherwise let $S:=S'\cup\{xw\}$. It is trivial to check that $S$ is a flippable set of $G$. Moreover, $|S|=|S'|+1\geq\third(n-4)+1>\third(n-2)$. Thus $G$ is not a counterexample.
\end{proof}

%%%%%%%%%%%%%%%%%%%%%%%%%%%%%%%%%%%%%%%%%%%%%%%%%%%%%%%%%%%%%%%%%%%%%%%%%%%%%%%

\begin{lemma} 
\lemlabel{FourFour}
A minimum counterexample has no edge $vw$ with $\deg(v)=4$ and $\deg(w)=4$.
\end{lemma}

\begin{proof} 
Let $G$ be a minimum counterexample with $n$ vertices. Suppose that $G$ has an edge $vw$ with $\deg(v)=4$ and $\deg(w)=4$. Let $b$ and $d$ be the vertices that see $vw$. Then $b\ne d$, as otherwise $G=K_3$. Let $a$ be the other neighbour of $v$. Let $c$ be the other neighbour of $w$. If $a=c$, then $G=K_5\setminus e$, in which case $G$ is not a counterexample by \lemref{SmallTriangBigFlip}. Thus $a\ne c$, and  $(a,b,c,d)$ is a $4$-cycle whose interior only contains $v$ and $w$. 

At least one of $ac$ and $bd$ is not an edge of $G$, as otherwise $G$ would contain a subdivision of $K_5$. If $ac$ is not an edge of $G$, then let $G'$ be the graph obtained from $G$ by deleting $v$ and $w$, and adding the edge $ac$. Otherwise, let $G'$ be the graph obtained from $G$ by deleting $v$ and $w$, and adding the edge $bd$. In both cases, $G'$ is a triangulation on $n-2$ vertices.
Since $G$ is minimum, $G'$ is not a counterexample. Thus $G'$ has a flippable set $S'$ of at least $\third(n-4)$ edges. Initialise $S:=S'$.

First suppose that $ac$ is not an edge of $G$. Then $ac$ is an edge of $G'$. 
If $ab\in S'$, then let $S:=S\cup\{wb\}$. 
If $bc\in S'$, then let $S:=S\cup\{vb\}$. 
If $cd\in S'$, then let $S:=S\cup\{vd\}$. 
If $ad\in S'$, then let $S:=S\cup\{wd\}$. 

Now suppose that $ac$ is an edge of $G$. Then $bd$ is an edge of $G'$. 
If $ab\in S'$, then let $S:=S'\cup\{vd\}$. 
If $ad\in S'$, then let $S:=S'\cup\{vb\}$. 
If $cd\in S'$, then let $S:=S'\cup\{wb\}$. 
If $bc\in S'$, then let $S:=S'\cup\{wd\}$. 

If none of these cases occur, then let  $S:=S\cup\{vb,wd\}$.
If both $vb$ and $vd$ have been added to $S$, then delete one from $S$.
If both $wb$ and $wd$ have been added to $S$, then delete one from $S$.
It is easily seen that in each case, $S$ is a flippable set, and
$|S|\geq|S'|+1\geq\third(n-4)+1>\third(n-2)$. Thus $G$ is not a counterexample.
\end{proof}

%%%%%%%%%%%%%%%%%%%%%%%%%%%%%%%%%%%%%%%%%%%%%%%%%%%%%%%%%%%%%%%%%%%%%%%%%%%%%%%

The following lemma is the key idea in the proof of \thmref{BigFlip}.

\begin{lemma}
\lemlabel{SiFlippable}
Let $\{E_1,E_2,E_3\}$ be an edge $3$-colouring of a triangulation $G$ such that every triangle is trichromatic. For each $1\leq i\leq 3$, let $S_i$ be the set of edges in $E_i$ that are not in a bad pair with some other edge in $E_i$. Then $S_i$ is flippable in $G$.
\end{lemma}

\begin{proof}
Since every triangle is trichromatic, no two edges in $S_i$ are consecutive. This is condition (1) in \lemref{FlippableCharacterisation}. Condition (2) in \lemref{FlippableCharacterisation} holds by the definition of $S_i$. Suppose that an edge $ab\in S_i$ is blocked by an edge $vw$. To show that condition (3) of \lemref{FlippableCharacterisation} is satisfied, we need to prove that $vw\in S_i$. 

First suppose that $vw\not\in E_i$. Since $(v,a,w)$ is a triangle, one of $av$ and $bv$ is in $E_i$, which implies that this edge and $ab$ are consecutive and both in $E_i$. This contradiction proves that $vw\in E_i$. Now suppose that $vw$ and some edge $xy$ form a bad pair. By \lemref{BlockingEdgeInBadPair}, $vw$ and $xy$ are in a common triangle. Thus $xy\not\in E_i$ and $vw$ does not form a bad pair with another edge in $E_i$. Therefore $vw\in S_i$, as desired. By \lemref{FlippableCharacterisation}, $S_i$ is flippable. 
\end{proof}

%%%%%%%%%%%%%%%%%%%%%%%%%%%%%%%%%%%%%%%%%%%%%%%%%%%%%%%%%%%%%%%%%%%%%%%%%%%%%

An edge is \emph{bad} if it is a member of a bad pair. An edge is \emph{good} if it is not bad. 

\begin{lemma}
\lemlabel{ReallyBadFace}
If every edge in a face $(u,v,w)$ of a triangulation $G$ is bad, then at least one of $\{u,v,w\}$ has degree three or four.
\end{lemma}

\begin{proof}
Assume $\deg(u)\leq\deg(v)\leq\deg(w)$. If $\deg(u)=3$ then we are done. Suppose that $\deg(u)\geq4$. Let $x,y,z$ be the other vertices that respectively see the edges $uv,vw,uw$. Since each of $u,v,w$ have degree at least four, $x,y,z$ are distinct. As illustrated in \figref{ReallyBadFace} with an outerface of $(u,v,w)$, there are edges $ab$, $cd$, and $ef$ such that $\{uv,ab\}$, $\{uw,ef\}$, $\{vw,cd\}$ are all bad pairs. For planarity to hold, and since $\deg(u)\leq\deg(v)\leq\deg(w)$, $d=x$ and $c=z$, which implies that $\deg(u)=4$, as desired.
\end{proof}

\Figure{ReallyBadFace}{\includegraphics{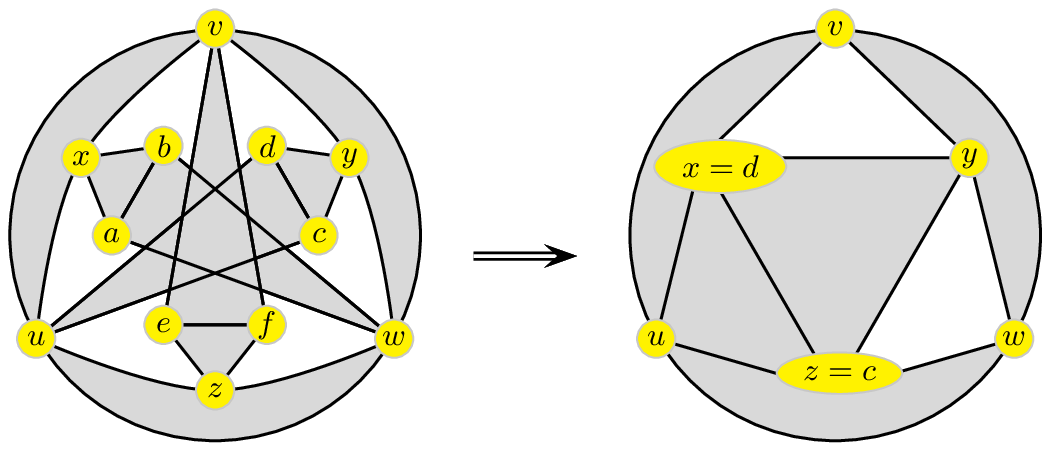}}{Three bad edges $uv$,  $uw$, and $vw$ all on one face.}

%%%%%%%%%%%%%%%%%%%%%%%%%%%%%%%%%%%%%%%%%%%%%%%%%%%%%%%%%%%%%%%%%%%%%%%%%%%%%%%%

\begin{lemma} 
\lemlabel{SeperatingCycle}
Define $S_1,S_2,S_3$ as in \lemref{SiFlippable}. Then every edge in a separating triangle is in $S_1\cup S_2\cup S_3$.
\end{lemma}

\begin{proof}
Consider an edge $vw\in E_i$ that is in a separating triangle $T$. If $vw$ is good then $vw\in S_i$ and we are done. Otherwise $vw$ is bad. By \lemref{EdgeInCycle}, the edge $e$ that forms a bad pair with $vw$ is also in $T$. Since each triangle is trichromatic, $e\not\in E_i$. Thus $vw\in S_i$. 
\end{proof}

%%%%%%%%%%%%%%%%%%%%%%%%%%%%%%%%%%%%%%%%%%%%%%%%%%%%%%%%%%%%%%%%%%%%%%%%%%%%%%%%

\begin{lemma} 
\lemlabel{BoundaryCycle}
In a minimum counterexample, every edge seen by a degree-$4$ vertex is
good.
\end{lemma}

\begin{proof} 
Let $v$ be a degree-$4$ vertex in a minimum counterexample $G$. 
Let $(a,b,c,d)$ be the neighbours of $v$ in cyclic order. 
Then $X:=\{ab,bc,cd,ad\}$ are the edges seen by $v$. 
Suppose on the contrary that one edge in $X$, say $ad$, is bad. 
Then $ad$ forms a bad pair with another edge in $X$.
Without loss of generality, either $\{ab,ad\}$ or $\{ad,bc\}$ are this bad pair.
If $\{ab,ad\}$ is a bad pair, then to avoid parallel edges, $\deg(a)=4$, which contradicts \lemref{FourFour}. 

Now suppose that $\{ad,bc\}$ is a bad pair. Let $x$ be the other vertex seen by these edges. Let $G'$ be the plane graph obtained from $G$  by deleting $v$, deleting the edges in the triangle $(c,d,x)$, merging the vertices $a$ and $d$, and merging the vertices $b$ and $c$. Then $G'$ is a triangulation on $n-3$ vertices. Since $G$ is minimum, $G'$ is not a counterexample. Thus $G'$ has a flippable set $S'$ of at least $\third(n-5)$ edges. Let $S:=S'\cup\{vd\}$. We claim that $S$ is flippable in $G$. Now $vd$ flips to $ac$, which is not an edge of $G$ as otherwise there would be a subdivision of $K_5$. The only edge that forms a bad pair with $vd$ is $vb$, which by construction is not in $S$. 
Thus $S$ is flippable, and $|S|=|S'|+1\geq\third(n-5)+1=\third(n-2)$. Thus $G$ is not a counterexample.
\end{proof}

%Now suppose that $\{ad,bc\}$ is a bad pair. Let $x$ be the other vertex seen by these edges. As illustrated in \figref{WeirdMerge}, let $G'$ be the plane graph obtained from $G$  by deleting $v$, deleting the edges in the triangle $(c,d,x)$,  merging the vertices $a$ and $c$, and merging the vertices $b$ and $d$. Then $G'$ is a triangulation on $n-1$ vertices. Since $G$ is minimum, $G'$ is not a counterexample. Thus $G'$ has a flippable set $S'$ of at least $\third(n-3)$ edges. 

%If $ab\in S'$, then let $S:=S'\cup\{cd\}$.  If $bx\in S'$, then let $S:=S'\cup\{dx\}$.  If $ax\in S'$, then let $S:=S'\cup\{cx\}$. If none of these cases occur, then let $S:=S'\cup\{vd\}$.

%It easily seen that in each case $S$ is a flippable set. Moreover, $|S|\geq|S'|+1\geq\third(n-3)+1>\third(n-2)$. Thus $G$ is not a counterexample.

%%%%%%%%%%%%%%%%%%%%%%%%%%%%%%%%%%%%%%%%%%%%%%%%%%%%%%%%%%%%%%%%%%%%%%%%%%

\begin{proof}[Proof of \thmref{BigFlip}]
Let $G$ be a minimum counterexample with $n$ vertices. By \lemref{TaitColouring}, there is a $3$-colouring $\{E_1,E_2,E_3\}$ of the edges of $G$ such that every triangle is trichromatic. Let $S_i$ be set of edges in $E_i$ that are not in a bad pair with another edge in $E_i$. By \lemref{SiFlippable}, $S_i$ is flippable. 

The neighbours of a degree-$3$ vertex form a separating triangle. By \lemref{SeperatingCycle}, every face incident to a degree-$3$ vertex has at least one edge in $S_1 \cup S_2 \cup S_3$. By \lemref{BoundaryCycle}, every face incident to a degree-$4$ vertex has one good edge, which is in $S_1 \cup S_2 \cup S_3$. By \lemref{ReallyBadFace}, every face not incident to degree-$3$ or degree-$4$ vertex has at least one good edge, which is in $S_1 \cup S_2 \cup S_3$. 

We conclude that every face has at least one edge in $S_1 \cup S_2 \cup S_3$. There are $2(n-2)$ faces and every edge is in two faces. Thus $|S_1 \cup S_2 \cup S_3|\geq n-2$. For some $i$, we have $|S_i| \geq\third(n-2)$. Therefore $G$ is not a counterexample, and since $G$ was minimum, there are no counterexamples.
\end{proof}

%%%%%%%%%%%%%%%%%%%%%%%%%%%%%%%%%%%%%%%%%%%%%%%%%%%%%%%%%%%%%%%%%%%%%%%%%%%%%
%\subsection{Upper Bounds}
%%%%%%%%%%%%%%%%%%%%%%%%%%%%%%%%%%%%%%%%%%%%%%%%%%%%%%%%%%%%%%%%%%%%%%%%%%%%%

Now for some upper bounds on \msf{G}.

\begin{lemma}
\lemlabel{BigFlipUniversalUpperBound}
For every $n$-vertex triangulation $G$, $\msf{G}\leq n-2$.
\end{lemma}

\begin{proof} Let $S$ be a flippable set of edges of $G$.  Every edge in $S$ is
incident to two distinct faces, and no other  edge on each of these faces is in
$S$.  (Otherwise there would be two consecutive edges in $S$.)\ There are
$2(n-2)$ faces in a triangulation. Thus $|S|\leq n-2$. \end{proof}

\begin{lemma}
\lemlabel{BigFlipExistentialUpperBound}  
There exist an $n$-vertex triangulation $G$ with $\msf{G}=\frac{6}{7}(n-2)$ for infinitely many $n$.
\end{lemma}

\begin{proof} 
Let $G_0$ be an arbitrary triangulation with $n_0$ vertices.  Let $G$ be the triangulation obtained from $G_0$ by adding a triangle inside each face $(u,v,w)$ of $G$, each vertex of which is adjacent to two of
$\{u,v,w\}$. Say $G$ has $n$ vertices. Then $n-2=n_0+3(2n_0-4)-2=7(n_0-2)$. Let $S$ be a flippable set of edges of $G$.  

There is at most one edge in $S$ on the boundary of each face of $G$. 
Suppose on the contrary that for some face $(u,v,w)$ of $G_0$, all seven of the 
corresponding faces of $G$ have an edge in $S$. Every edge in $S$ is on the boundary of two faces of $G$. Thus $|S\cap\{uv,uw,vw\}|=1$ or $3$. Let $(x,y,z)$ be the triangle of $G$ inside $(u,v,w)$, with connecting edges $\{xv,xw,yu,yw,zu,zv\}$.

Case 1. $|S\cap\{uv,uw,vw\}|=1$: Without loss of generality $S\cap\{uv,uw,vw\}=\{uv\}$, as illustrated in \figref{UpperBound}(a) and (b). Thus either (a) $uy\in S$ or (b) $zy\in S$. If $uy\in S$, then $xy\not\in S$ (as otherwise \F{G}{S} would have parallel edges). Thus $\{xz,xw\}\in S$, in which case $\F{G}{S}$ has parallel edges, a contradiction. If $zy\in S$, then $yw\in S$, as otherwise no edge on $(u,w,y)$ would be in $S$. In this case \F{G}{S} has parallel edges.

Case 2. $|S\cap\{uv,uw,vw\}|=3$: Then $zy$ is the only edge on the boundary of the face $(u,z,y)$ that can be flipped, as illustrated in \figref{UpperBound}(c). Hence $zy\in S$. This implies that no edge on the faces $(z,v,x)$ and $(x,y,w)$ can be flipped, a contradiction.

\Figure{UpperBound}{\includegraphics{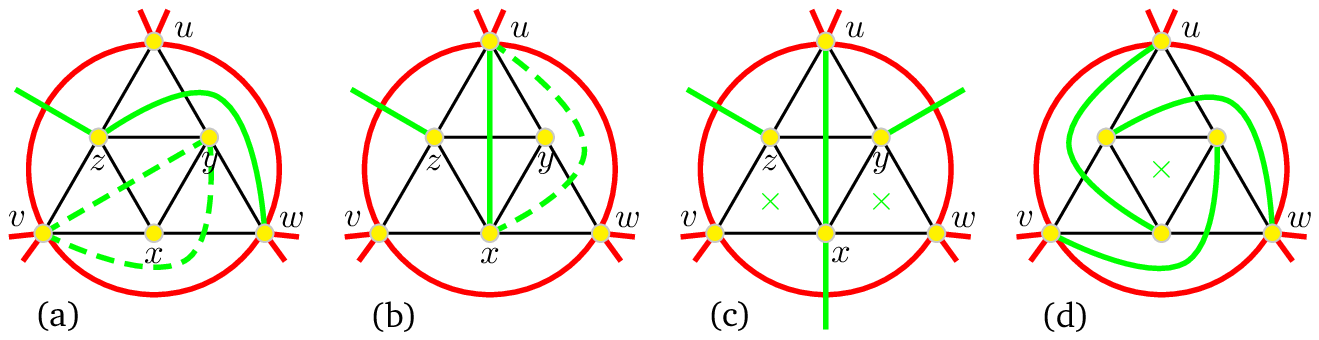}}{(a)--(c) For any number of flips in the outer triangle, at least one internal face does not have an edge in $S$. (d) How to construct a flip set for $G$.}

Therefore for every face of $G_0$, at least one of the seven corresponding faces of $G$ does not have an edge in $S$. Hence at least $2(n_0-2)=\frac{2}{7}(n-2)$ faces of $G$ do not have an edge in $S$.  Every face of $G$ has at most one edge in $S$. Thus $|S|\leq\half(2(n-2)-\frac{2}{7}(n-2))=\frac{6}{7}(n-2)$. 

It remains to construct a flippable set of $\frac{6}{7}(n-2)$ edges in $G$. 
For each face of $G_0$, add the edges shown in \figref{UpperBound}(d)  to a
set $S$.  Clearly $S$ is flippable. In every face of $G_0$, exactly one of the
corresponding seven faces of $G$ does not have an edge in $S$, and the
remaining six faces each have exactly one edge in $S$. By the above analysis,
$|S|=\frac{6}{7}(n-2)$. \end{proof}

An obvious open problem is to close the gap between the lower bound of 
$\third(n-2)$ and the upper bound of $\frac{6}{7}(n-2)$ in the above
results. For $5$-connected triangulations we can improve the lower bound as follows.

\begin{theorem}
\thmlabel{BigFlipFiveConnLowerBound}
For every $5$-connected triangulation $G$ with $n$ vertices, $\msf{G}=n-2$.
\end{theorem}

\begin{proof} 
Observe that every edge in $G$ is flippable, as otherwise $G$ has a separating triangle (since $G$ has at least five vertices). There is no bad pair in $G$, as otherwise $G$ has a separating $4$-cycle. By \lemref{FlippableCharacterisation}, a set of edges $S$ in a $5$-connected triangulation $G$ is flippable if and only if no two edges in $S$ are consecutive. By  \lemref{TriangleCover}, $G$ has a set of edges $S$ such that every triangle of $G$ has exactly one edge in $S$. Thus no two edges in $S$ are consecutive. Hence $S$ is flippable. By the argument employed in  \lemref{BigFlipUniversalUpperBound}, $|S|=n-2$.  Therefore $\msf{G}\geq n-2$. By \lemref{BigFlipUniversalUpperBound}, $\msf{G}\leq n-2$. \end{proof}

%%%%%%%%%%%%%%%%%%%%%%%%%%%%%%%%%%%%%%%%%%%%%%%%%%%%%%%%%%%%%%%%%%%%%%%%%%%%%%%
\subsubsection*{Acknowledgements} 

The research of Prosenjit Bose was partially completed at the Departament de Matem{\`a}tica Aplicada II, Universitat Polit{\`e}cnica de Catalunya, Barcelona, Spain. Thanks to Ferran Hurtado for graciously hosting the visit. Thanks to the referees for pointing out an error in a preliminary version of the paper. 

%%%%%%%%%%%%%%%%%%%%%%%%%%%%%%%%%%%%%%%%%%%%%%%%%%%%%%%%%%%%%%%%%%%%%%%%%%%%%%%%
%\bibliographystyle{myNatbibStyle}
%\bibliography{myBibliography,myConferences}

\begin{thebibliography}{37}
\providecommand{\natexlab}[1]{#1}
\providecommand{\url}[1]{\texttt{#1}}
\providecommand{\urlprefix}{}
\expandafter\ifx\csname urlstyle\endcsname\relax
  \providecommand{\doi}[1]{doi:\discretionary{}{}{}#1}\else
  \providecommand{\doi}{doi:\discretionary{}{}{}\begingroup
  \urlstyle{rm}\Url}\fi

\bibitem[{Biedl et~al.(2001)Biedl, Bose, Demaine, and Lubiw}]{BBDL-JAlg01}
\textsc{Therese~C. Biedl, Prosenjit Bose, Erik~D. Demaine, and Anna Lubiw}.
\newblock Efficient algorithms for {P}etersen's matching theorem.
\newblock \emph{J. Algorithms}, 38(1):110--134, 2001.

\bibitem[{Bose et~al.(2006)Bose, Dujmovi\'c, and Wood}]{BDW-CDM06}
\textsc{Prosenjit Bose, Vida Dujmovi\'c, and David~R. Wood}.
\newblock Induced subgraphs of bounded treewidth and bounded degree.
\newblock \emph{Contrib. Discrete Math.}, 1(1):88--105, 2006.

\bibitem[{Brunet et~al.(1996)Brunet, Nakamoto, and Negami}]{BNN-JCTB96}
\textsc{Richard Brunet, Atsuhiro Nakamoto, and Seiya Negami}.
\newblock Diagonal flips of triangulations on closed surfaces preserving
  specified properties.
\newblock \emph{J. Combin. Theory Ser. B}, 68(2):295--309, 1996.

\bibitem[{Chen et~al.(2005)Chen, Chang, and Wang}]{CCW-IJCM05}
\textsc{Yen-Ju Chen, Jou-Ming Chang, and Yue-Li Wang}.
\newblock An efficient algorithm for estimating rotation distance between two
  binary trees.
\newblock \emph{Int. J. Comput. Math.}, 82(9):1095--1106, 2005.

\bibitem[{Chiba and Nishizeki(1989)}]{CN89}
\textsc{Norishige Chiba and Takao Nishizeki}.
\newblock The {H}amiltonian cycle problem is linear-time solvable for
  {$4$}-connected planar graphs.
\newblock \emph{J. Algorithms}, 10(2):187--211, 1989.

\bibitem[{Chrobak and Payne(1995)}]{ChrobakPayne-IPL95}
\textsc{Marek Chrobak and Thomas~H. Payne}.
\newblock A linear-time algorithm for drawing a planar graph on a grid.
\newblock \emph{Inform. Process. Lett.}, 54(4):241--246, 1995.

\bibitem[{Cort\'{e}s et~al.(2002)Cort\'{e}s, Grima, Marquez, and
  Nakamoto}]{CGMN-DM02}
\textsc{Carmen Cort\'{e}s, Clara Grima, Alberto Marquez, and Atsuhiro
  Nakamoto}.
\newblock Diagonal flips in outer-triangulations on closed surfaces.
\newblock \emph{Discrete Math.}, 254(1-3):63--74, 2002.

\bibitem[{Cort\'{e}s and Nakamoto(2000)}]{CN-DM00a}
\textsc{Carmen Cort\'{e}s and Atsuhiro Nakamoto}.
\newblock Diagonal flips in outer-torus triangulations.
\newblock \emph{Discrete Math.}, 216(1-3):71--83, 2000.

\bibitem[{{de Fraysseix} et~al.(1990){de Fraysseix}, Pach, and
  Pollack}]{dFPP90}
\textsc{Hubert {de Fraysseix}, J\'{a}nos Pach, and Richard Pollack}.
\newblock How to draw a planar graph on a grid.
\newblock \emph{Combinatorica}, 10(1):41--51, 1990.

\bibitem[{Galtier et~al.(2003)Galtier, Hurtado, Noy, P\'{e}rennes, and
  Urrutia}]{GHNP-IJCGA03}
\textsc{Jer{\^o}me Galtier, Ferran Hurtado, Marc Noy, Stephane P\'{e}rennes,
  and Jorge Urrutia}.
\newblock Simultaneous edge flipping in triangulations.
\newblock \emph{Internat. J. Comput. Geom. Appl.}, 13(2):113--133, 2003.

\bibitem[{Gao et~al.(2001)Gao, Urrutia, and Wang}]{GUW-GC01}
\textsc{Zhicheng Gao, Jorge Urrutia, and Jianyu Wang}.
\newblock Diagonal flips in labelled planar triangulations.
\newblock \emph{Graphs Combin.}, 17(4):647--657, 2001.

\bibitem[{Gao and Wang(1999)}]{GW-JCTA99}
\textsc{Zhicheng Gao and Jianyu Wang}.
\newblock Enumeration of rooted planar triangulations with respect to diagonal
  flips.
\newblock \emph{J. Combin. Theory Ser. A}, 88(2):276--296, 1999.

\bibitem[{Hanke et~al.(1996)Hanke, Ottmann, and Schuierer}]{HOS-JUCS96}
\textsc{Sabine Hanke, Thomas Ottmann, and Sven Schuierer}.
\newblock The edge-flipping distance of triangulations.
\newblock \emph{J.UCS}, 2(8):570--579, 1996.

\bibitem[{Hurtado and Noy(1999)}]{HN-CGTA99}
\textsc{Ferran Hurtado and Marc Noy}.
\newblock Graph of triangulations of a convex polygon and tree of
  triangulations.
\newblock \emph{Comput. Geom.}, 13(3):179--188, 1999.

\bibitem[{Hurtado et~al.(1999)Hurtado, Noy, and Urrutia}]{HNU-DCG99}
\textsc{Ferran Hurtado, Marc Noy, and Jorge Urrutia}.
\newblock Flipping edges in triangulations.
\newblock \emph{Discrete Comput. Geom.}, 22(3):333--346, 1999.

\bibitem[{Komuro(1997)}]{Komuro-Yoko97}
\textsc{Hideo Komuro}.
\newblock The diagonal flips of triangulations on the sphere.
\newblock \emph{Yokohama Math. J.}, 44(2):115--122, 1997.

\bibitem[{Komuro and Ando(2001)}]{KA-AC01}
\textsc{Hideo Komuro and Kiyoshi Ando}.
\newblock Diagonal flips of pseudo triangulations on the sphere.
\newblock \emph{Ars Combin.}, 59:225--239, 2001.

\bibitem[{Komuro et~al.(1999)Komuro, Nakamoto, and Negami}]{KNN-JCTB99}
\textsc{Hideo Komuro, Atsuhiro Nakamoto, and Seiya Negami}.
\newblock Diagonal flips in triangulations on closed surfaces with minimum
  degree at least {$4$}.
\newblock \emph{J. Combin. Theory Ser. B}, 76(1):68--92, 1999.

\bibitem[{K{\"o}nig(1936)}]{Konig36}
\textsc{D{\'e}nes K{\"o}nig}.
\newblock \emph{Theorie der endlichen und unendlichen {G}raphen.
  {K}ombinatorische {T}opologie der {S}treckenkomplexe}.
\newblock Akademische Verlagsgesellschaft, Leipzig, 1936.

\bibitem[{Mori et~al.(2003)Mori, Nakamoto, and Ota}]{MNO-GC03}
\textsc{Ryuichi Mori, Atsuhiro Nakamoto, and Katsuhiro Ota}.
\newblock Diagonal flips in {H}amiltonian triangulations on the sphere.
\newblock \emph{Graphs Combin.}, 19(3):413--418, 2003.

\bibitem[{Nakamigawa(2000)}]{Nakamigawa-TCS00}
\textsc{Tomoki Nakamigawa}.
\newblock A generalization of diagonal flips in a convex polygon.
\newblock \emph{Theoret. Comput. Sci.}, 235(2):271--282, 2000.

\bibitem[{Nakamoto and Negami(2002)}]{NN-Yoko02}
\textsc{Atsuhiro Nakamoto and Seiya Negami}.
\newblock Diagonal flips in graphs on closed surfaces with specified face size
  distributions.
\newblock \emph{Yokohama Math. J.}, 49(2):171--180, 2002.

\bibitem[{Nakamoto et~al.(2006)Nakamoto, Sakuma, and Suzuki}]{NSS-JGT06}
\textsc{Atsuhiro Nakamoto, Tadashi Sakuma, and Yusuke Suzuki}.
\newblock {$N$}-flips in even triangulations on the sphere.
\newblock \emph{J. Graph Theory}, 51(3):260--268, 2006.

\bibitem[{Negami(1994)}]{Negami-DM94}
\textsc{Seiya Negami}.
\newblock Diagonal flips in triangulations of surfaces.
\newblock \emph{Discrete Math.}, 135(1-3):225--232, 1994.

\bibitem[{Negami(1998)}]{Negami-Yoko98}
\textsc{Seiya Negami}.
\newblock Diagonal flips in triangulations on closed surfaces, estimating upper
  bounds.
\newblock \emph{Yokohama Math. J.}, 45(2):113--124, 1998.

\bibitem[{Negami(1999)}]{Negami99}
\textsc{Seiya Negami}.
\newblock Diagonal flips of triangulations on surfaces, a survey.
\newblock \emph{Yokohama Math. J.}, 47:1--40, 1999.

\bibitem[{Negami and Nakamoto(1993)}]{NN93}
\textsc{Seiya Negami and Atsuhiro Nakamoto}.
\newblock Diagonal transformations of graphs on closed surfaces.
\newblock \emph{Sci. Rep. Yokohama Nat. Univ. Sect. I Math. Phys. Chem.},
  40:71--97, 1993.

\bibitem[{Pallo(1987)}]{Pallo-IPL87}
\textsc{Jean Pallo}.
\newblock On the rotation distance in the lattice of binary trees.
\newblock \emph{Inform. Process. Lett.}, 25(6):369--373, 1987.

\bibitem[{Pallo(2000)}]{Pallo-IPL00}
\textsc{Jean Pallo}.
\newblock An efficient upper bound of the rotation distance of binary trees.
\newblock \emph{Inform. Process. Lett.}, 73(3-4):87--92, 2000.

\bibitem[{Petersen(1891)}]{Petersen1891}
\textsc{Julius Petersen}.
\newblock {D}ie {T}heorie der regul{\"a}ren {G}raphen.
\newblock \emph{Acta. Math.}, 15:193--220, 1891.

\bibitem[{Robertson et~al.(1997)Robertson, Sanders, Seymour, and
  Thomas}]{RSST97}
\textsc{Neil Robertson, Daniel~P. Sanders, Paul~D. Seymour, and Robin Thomas}.
\newblock The four-colour theorem.
\newblock \emph{J. Combin. Theory Ser. B}, 70(1):2--44, 1997.

\bibitem[{Sleator et~al.(1988)Sleator, Tarjan, and Thurston}]{STT-JAMS88}
\textsc{Daniel~D. Sleator, Robert~E. Tarjan, and William~P. Thurston}.
\newblock Rotation distance, triangulations, and hyperbolic geometry.
\newblock \emph{J. Amer. Math. Soc.}, 1(3):647--681, 1988.

\bibitem[{Sleator et~al.(1992)Sleator, Tarjan, and Thurston}]{STT-SJDM92}
\textsc{Daniel~D. Sleator, Robert~E. Tarjan, and William~P. Thurston}.
\newblock Short encodings of evolving structures.
\newblock \emph{SIAM J. Discrete Math.}, 5(3):428--450, 1992.

\bibitem[{Tait(1880)}]{Tait1880a}
\textsc{Peter~Guthrie Tait}.
\newblock Note on a theorem in geometry of position.
\newblock \emph{Trans. Roy. Soc. Edinburgh}, 29:657--660, 1880.

\bibitem[{Wagner(1936)}]{Wagner36}
\textsc{Klaus Wagner}.
\newblock Bemerkung zum {V}ierfarbenproblem.
\newblock \emph{Jber. Deutsch. Math.-Verein.}, 46:26--32, 1936.

\bibitem[{Watanabe and Negami(1999)}]{WN-Yoko99}
\textsc{Takahiro Watanabe and Seiya Negami}.
\newblock Diagonal flips in pseudo-triangulations on closed surfaces without
  loops.
\newblock \emph{Yokohama Math. J.}, 47:213--223, 1999.

\bibitem[{Whitney(1931)}]{Whitney31}
\textsc{Hassler Whitney}.
\newblock A theorem on graphs.
\newblock \emph{Ann. of Math. (2)}, 32(2):378--390, 1931.

\end{thebibliography}
%%%%%%%%%%%%%%%%%%%%%%%%%%%%%%%%%%%%%%%%%%%%%%%%%%%%%%%%%%%%%%%%%%%%%%%%%%%%%%%%

\def\soft#1{\leavevmode\setbox0=\hbox{h}\dimen7=\ht0\advance \dimen7
  by-1ex\relax\if t#1\relax\rlap{\raise.6\dimen7
  \hbox{\kern.3ex\char'47}}#1\relax\else\if T#1\relax
  \rlap{\raise.5\dimen7\hbox{\kern1.3ex\char'47}}#1\relax \else\if
  d#1\relax\rlap{\raise.5\dimen7\hbox{\kern.9ex \char'47}}#1\relax\else\if
  D#1\relax\rlap{\raise.5\dimen7 \hbox{\kern1.4ex\char'47}}#1\relax\else\if
  l#1\relax \rlap{\raise.5\dimen7\hbox{\kern.4ex\char'47}}#1\relax \else\if
  L#1\relax\rlap{\raise.5\dimen7\hbox{\kern.7ex
  \char'47}}#1\relax\else\message{accent \string\soft \space #1 not
  defined!}#1\relax\fi\fi\fi\fi\fi\fi} \def\cprime{$'$}

\end{document}